\documentstyle[12pt,amssymb,amscd]{article}
\def\frak{{\cal}}

{
\begin{equation}}
{\end{equation}%
\addtocounter{equation}{-1}}

\begin{document}

\baselineskip 18pt

\newtheorem{definition}{\normalsize\sc Definition}[section]
\newtheorem{defin}[definition]{\normalsize\sc Definition}
\newtheorem{prop}[definition]{\normalsize\sc Proposition}
\newtheorem{lem}[definition]{\normalsize\sc Lemma}
\newtheorem{corollary}[definition]{\normalsize\sc Corollary}
\newtheorem{corol}[definition]{\normalsize\sc Corollary}
\newtheorem{theorem}[definition]{\normalsize\sc Theorem}
\newtheorem{example}[definition]{\normalsize\sc  Example}
\newtheorem{remark}[definition]{\normalsize\sc Remark}

\newcommand{\nc}[2]{\newcommand{#1}{#2}}
\newcommand{\rnc}[2]{\renewcommand{#1}{#2}}

\nc{\Section}{\setcounter{definition}{0}\section}
\newcounter{c}
\renewcommand{\[}{\setcounter{c}{1}$$}
\newcommand{\etyk}[1]{\vspace{-7.4mm}$$\begin{equation}\Label{#1}
\addtocounter{c}{1}}
\renewcommand{\]}{\ifnum \value{c}=1 $$\else \end{equation}\fi}

\nc{\bpr}{\begin{prop}}
\nc{\bth}{\begin{theorem}}
\nc{\ble}{\begin{lem}}
\nc{\bco}{\begin{corollary}}
\nc{\bre}{\begin{remark}}
\nc{\bex}{\begin{example}}
\nc{\bde}{\begin{definition}}
\nc{\ede}{\end{definition}}
\nc{\epr}{\end{prop}}
\nc{\ethe}{\end{theorem}}
\nc{\ele}{\end{lem}}
\nc{\eco}{\end{corollary}}
\nc{\ere}{\hfill\mbox{$\Diamond$}\end{remark}}
\nc{\eex}{\end{example}}
\nc{\epf}{\hfill\mbox{$\Box$}}
\nc{\ot}{\otimes}
\nc{\bsb}{\begin{Sb}}
\nc{\esb}{\end{Sb}}
\nc{\ct}{\mbox{${\cal T}$}}
\nc{\ctb}{\mbox{${\cal T}\sb B$}}
\nc{\bcd}{\[\begin{CD}}
\nc{\ecd}{\end{CD}\]}
\nc{\ba}{\begin{array}}
\nc{\ea}{\end{array}}
\nc{\bea}{\begin{eqnarray}}
\nc{\eea}{\end{eqnarray}}
\nc{\be}{\begin{enumerate}}
\nc{\ee}{\end{enumerate}}
\nc{\beq}{\begin{equation}}
\nc{\eeq}{\end{equation}}
\nc{\bi}{\begin{itemize}}
\nc{\ei}{\end{itemize}}
\nc{\kr}{\mbox{Ker}}
\nc{\te}{\!\ot\!}
\nc{\pf}{\mbox{$P\!\sb F$}}
\nc{\pn}{\mbox{$P\!\sb\nu$}}
\nc{\bmlp}{\mbox{\boldmath$\left(\right.$}}
\nc{\bmrp}{\mbox{\boldmath$\left.\right)$}}
\rnc{\phi}{\mbox{$\varphi$}}
\nc{\LAblp}{\mbox{\LARGE\boldmath$($}}
\nc{\LAbrp}{\mbox{\LARGE\boldmath$)$}}
\nc{\Lblp}{\mbox{\Large\boldmath$($}}
\nc{\Lbrp}{\mbox{\Large\boldmath$)$}}
\nc{\lblp}{\mbox{\large\boldmath$($}}
\nc{\lbrp}{\mbox{\large\boldmath$)$}}
\nc{\blp}{\mbox{\boldmath$($}}
\nc{\brp}{\mbox{\boldmath$)$}}
\nc{\LAlp}{\mbox{\LARGE $($}}
\nc{\LArp}{\mbox{\LARGE $)$}}
\nc{\Llp}{\mbox{\Large $($}}
\nc{\Lrp}{\mbox{\Large $)$}}
\nc{\llp}{\mbox{\large $($}}
\nc{\lrp}{\mbox{\large $)$}}
\nc{\lbc}{\mbox{\Large\boldmath$,$}}
\nc{\lc}{\mbox{\Large$,$}}
\nc{\Lall}{\mbox{\Large$\forall$}}
\nc{\bc}{\mbox{\boldmath$,$}}
\rnc{\epsilon}{\varepsilon}
\nc{\ra}{\rightarrow}
\nc{\ci}{\circ}
\nc{\cc}{\!\ci\!}
\nc{\T}{\mbox{\sf T}}
\nc{\can}{\mbox{\em\sf T}\!\sb R}
\nc{\cnl}{$\mbox{\sf T}\!\sb R$}
\nc{\lra}{\longrightarrow}
\nc{\M}{\mbox{Map}}
\nc{\imp}{\Rightarrow}
\rnc{\iff}{\Leftrightarrow}
\nc{\bmq}{\cite{bmq}}
\nc{\ob}{\mbox{$\Omega\sp{1}\! (\! B)$}}
\nc{\op}{\mbox{$\Omega\sp{1}\! (\! P)$}}
\nc{\oa}{\mbox{$\Omega\sp{1}\! (\! A)$}}
\nc{\inc}{\mbox{$\,\subseteq\;$}}
\nc{\de}{\mbox{$\Delta$}}
\nc{\spp}{\mbox{${\cal S}{\cal P}(P)$}}
\nc{\dr}{\mbox{$\Delta_{R}$}}
\nc{\dsr}{\mbox{$\Delta_{\cal R}$}}
\nc{\m}{\mbox{m}}
\nc{\0}{\sb{(0)}}
\nc{\1}{\sb{(1)}}
\nc{\2}{\sb{(2)}}
\nc{\3}{\sb{(3)}}
\nc{\4}{\sb{(4)}}
\nc{\5}{\sb{(5)}}
\nc{\6}{\sb{(6)}}
\nc{\7}{\sb{(7)}}
\nc{\hsp}{\hspace*}
\nc{\nin}{\mbox{$n\in\{ 0\}\!\cup\!{\Bbb N}$}}
\nc{\al}{\mbox{$\alpha$}}
\nc{\bet}{\mbox{$\beta$}}
\nc{\ha}{\mbox{$\alpha$}}
\nc{\hb}{\mbox{$\beta$}}
\nc{\hg}{\mbox{$\gamma$}}
\nc{\hd}{\mbox{$\delta$}}
\nc{\he}{\mbox{$\varepsilon$}}
\nc{\hz}{\mbox{$\zeta$}}
\nc{\hs}{\mbox{$\sigma$}}
\nc{\hk}{\mbox{$\kappa$}}
\nc{\hm}{\mbox{$\mu$}}
\nc{\hn}{\mbox{$\nu$}}
\nc{\la}{\mbox{$\lambda$}}
\nc{\hl}{\mbox{$\lambda$}}
\nc{\hG}{\mbox{$\Gamma$}}
\nc{\hD}{\mbox{$\Delta$}}
\nc{\ho}{\mbox{$\omega$}}
\nc{\hO}{\mbox{$\Omega$}}
\nc{\hp}{\mbox{$\pi$}}
\nc{\hP}{\mbox{$\Pi$}}
\nc{\bpf}{{\it Proof.~~}}
\nc{\slq}{\mbox{$A(SL\sb q(2))$}}
\nc{\fr}{\mbox{$Fr\llp A(SL(2,\IC))\lrp$}}
\nc{\slc}{\mbox{$A(SL(2,\IC))$}}
\nc{\af}{\mbox{$A(F)$}}
\rnc{\widetilde}{\tilde}
\nc{\qdt}{quantum double torus}
\nc{\aqdt}{\mbox{$A(DT^2_q)$}}
\nc{\dtq}{\mbox{$DT^2_q$}}
\nc{\uc}{\mbox{$U(2)$}}
\nc{\uq}{\mbox{$U_{q^{-1},q}(2)$}}
\rnc{\subset}{\inc}

\def\esl{{\mbox{$E\sb{\frak s\frak l (2,{\Bbb C})}$}}}
\def\esu{{\mbox{$E\sb{\frak s\frak u(2)}$}}}
\def\zf{{\mbox{${\Bbb Z}\sb 4$}}}
\def\zt{{\mbox{$2{\Bbb Z}\sb 2$}}}
\def\ox{{\mbox{$\Omega\sp 1\sb{\frak M}X$}}}
\def\oxh{{\mbox{$\Omega\sp 1\sb{\frak M-hor}X$}}}
\def\oxs{{\mbox{$\Omega\sp 1\sb{\frak M-shor}X$}}}
\def\Fr{\mbox{Fr}}
\def\gal{-Galois extension}
\def\hge{Hopf-Galois extension}
\def\cge{coalgebra-Galois extension}
\def\pge{$\psi$-Galois extension}
\def\ta{\tilde a}
\def\tb{\tilde b}
\def\tc{\tilde c}
\def\td{\tilde d}
\def\st{\stackrel}

\newcommand{\Sp}{{\rm Sp}\,}
\newcommand{\Mor}{\mbox{$\rm Mor$}}
\newcommand{\skrA}{{\cal A}}
\newcommand{\Phase}{\mbox{$\rm Phase\,$}}
\newcommand{\id}{{\rm id}}
\newcommand{\diag}{{\rm diag}}
\newcommand{\inv}{{\rm inv}}
\newcommand{\ad}{{\rm ad}}
\newcommand{\poi}{{\rm pt}}
\newcommand{\Dim}{{\rm dim}\,}
\newcommand{\Ker}{{\rm ker}\,}
\newcommand{\Mat}{{\rm Mat}\,}
\newcommand{\Rep}{{\rm Rep}\,}
\newcommand{\Fun}{{\rm Fun}\,}
\newcommand{\Tr}{{\rm Tr}\,}
\newcommand{\supp}{\mbox{$\rm supp$}}
\newcommand{\half}{\frac{1}{2}}
\newcommand{\skrF}{{A}}
\newcommand{\skrD}{{\cal D}}
\newcommand{\skrC}{{\cal C}}
\newcommand{\ttimes}{\mbox{$\hspace{.5mm}\bigcirc\hspace{-4.9mm}
\perp\hspace{1mm}$}}
\newcommand{\Ttimes}{\mbox{$\hspace{.5mm}\bigcirc\hspace{-3.7mm}
\raisebox{-.7mm}{$\top$}\hspace{1mm}$}}
\newcommand{\bbr}{{\bf R}}
\newcommand{\bbz}{{\bf Z}}
\newcommand{\Ci}{C_{\infty}}
\newcommand{\Cb}{C_{b}}
\newcommand{\fa}{\forall}
\newcommand{\rrr}{right regular representation}
\newcommand{\wrt}{with respect to}
\newcommand{\qg}{quantum group}
\newcommand{\qgs}{quantum groups}
\newcommand{\cs}{classical space}
\newcommand{\qs}{quantum space}
\newcommand{\po}{Pontryagin}
\newcommand{\ch}{character}
\newcommand{\chs}{characters}

\def\inbar{\,\vrule height1.5ex width.4pt depth0pt}
\def\IC{{\Bbb C}}
\def\IZ{{\Bbb Z}}
\def\IN{{\Bbb N}}
\def\otc{\otimes_{\IC}}
\def\ra{\rightarrow}
\def\ota{\otimes_ A}
\def\otza{\otimes_{ Z(A)}}
\def\otc{\otimes_{\IC}}
\def\h{\rho}
\def\x{\zeta}
\def\th{\theta}
\def\s{\sigma}
\def\t{\tau}
\def\te{{\tilde e}}

\def\sw#1{{\sb{(#1)}}}
\def\sco#1{{\sp{(\bar #1)}}} 
\def\su#1{{\sp{(#1)}}} 
\def\extd{{\rm d}}
\def\covD{{\rm D}}
\def\proof{\noindent{\sl Proof.~~}}
\def\endproof{\hbox{$\sqcup$}\llap{\hbox{$\sqcap$}}\medskip}
\def\tens{\mathop{\otimes}}
\def\CC{{\frak C}}
\def\CL{{\Lambda}}
\def\M{{\bf M}}
\def\Mod{{\rm  Mod}}
\def\o{{}_{(1)}}
\def\t{{}_{(2)}}
\def\Bo{{}_{\und{(1)}}}
\def\Bt{{}_{\und{(2)}}}
\def\th{{}_{(3)}}
\def\Bth{{}_{\und{(3)}}}
\def\<{{\langle}}
\def\>{{\rangle}}
\def\und#1{{\underline{#1}}}
\def\ra{{\triangleleft}}
\def\la{{\triangleright}} 
\def\id{{\rm id}} 
\def\eps{\epsilon}
\def\span{{\rm span}}
\def\q2{{q^{-2}}}
\def\bicross{{\blacktriangleright\!\!\!\triangleleft}}
\def\note#1{{}}
\def\eqn#1#2{\begin{equation}#2\label{#1}\end{equation}}
\def\qbinom#1#2#3{\left(\begin{array}{c}#1\\#2\end{array}\right)\sb#3}
\def\Z{{\Bbb Z}}
\def\can{{\rm can}}
\def\note#1{}
\def\CA{{\frak A}}

\def\align#1{\begin{eqnarray*}#1\end{eqnarray*}}
\def\cmath#1{\[\begin{array}{c} #1 \end{array}\]}
\def\ceqn#1#2{\begin{equation}\label{#1}
\begin{array}{c}#2\end{array}\end{equation}}
\def\equad{\kern -1.7em}
\def\qqquad{\qquad\quad}
\def\nquad{{\!\!\!\!\!\!}}
\def\nqquad{\nquad\nquad}\def\nqqquad{\nqquad\nquad}
\def\C{{\Bbb C}}\def\H{{\Bbb H}}
\def\uo{{}^{(1)}}\def\uth{{}^{(3)}}
\def\Ad{{\rm Ad}}
\def\Hom{{\rm Hom}}
\def\isom{{\cong}}
\def\trace{{\rm Tr}}
\def\ut{{}^{(2)}}\def\umo{{}^{-(1)}}\def\umt{{}^{-(2)}}
\def\cg{{\sl g}}
\def\R{{\Bbb R}}
\def\tr{{\rm Tr}}
\def\CM{{\cal M}}
\def\CN{{\cal N}}
\def\CQ{{\cal Q}}

\begin{center}
{\large\bf QUANTUM GEOMETRY OF ALGEBRA FACTORISATIONS AND COALGEBRA
BUNDLES}
\vspace{10pt}\\
Tomasz Brzezi\'nski\footnote{EPSRC Advanced Research Fellow at the
University of Wales Swansea.}\vspace{6pt}\\
\normalsize Department of Mathematics, University of Wales Swansea,\\
\normalsize Singleton Park, Swansea SA2 8PP, U.K.\\
and\\
Department of Theoretical Physics, University of \L\'od\'z,\\
Pomorska 149/153, 90--236 \L\'od\'z, Poland
\vspace{6pt}\\
+\vspace{6pt}\\
Shahn Majid\footnote{Reader and Royal Society University Research
Fellow at QMW and Senior Research Fellow at the Department of Applied
Mathematics and Theoretical Physics, University of Cambridge, England.}
\vspace{6pt}\\
\normalsize School of Mathematical Sciences, Queen Mary and Westfield
College,\\ Mile End Road, London E1 4NS, England.
\end{center}

\begin{center}
\end{center}

\begin{abstract} We develop the noncommutative geometry (bundles, 
connections etc.)
associated to algebras that factorise into two subalgebras. 
An example is the factorisation of matrices $M_2(\C)=\C\Z_2\cdot\C\Z_2$.
We also further extend the coalgebra version of theory introduced
previously, to include frame resolutions and corresponding covariant
derivatives and torsions.  As an example, we construct $q$-monopoles on
all the  Podle\'s quantum spheres $S^2_{q,s}$.

\end{abstract}

\section*{\normalsize\sc\centering 1. Introduction}\vspace{-.6\baselineskip}
\setcounter{section}{1}

In \cite{BrzMa:coa} it was shown that one can generalise the notion
of principal bundles in noncommutative geometry\cite{BrzMa:gau} to
a very general setting in which the role of   `coordinate
functions' on the base is played by a general (possibly
noncommutative) algebra and the role of the `structure group'
(fibre) of the principal bundle is a coalgebra. In particular, it
need not be a quantum group, which would be too restrictive for
many interesting examples. In \cite{Brz:mod} the theory of modules
or `associated bundles' is extended to this case along the lines of
the quantum group case in \cite{BrzMa:gau}. We apply this now to
extend the recently introduced notion of a frame resolution
\cite{Ma:rie}, thereby bringing the coalgebra version of the gauge theory
in line with the more restrictive quantum group gauge theory case.

The paper begins, however, in Section~2 with a useful reformulation of
coalgebra bundles entirely in terms of algebras. This is a  theory where the
role of `gauge group' or fibre in the principal bundle is played by any algebra 
$A$ subject to a certain nondegeneracy  `Galois' condition for its action on the 
algebra $P$ of the total space of the
bundle. The algebra $A$ plays the role of a classical (or quantum) enveloping algebra 
of a Lie algebra in usual (or quantum group) gauge theory, but now without any kind 
of Hopf algebra structure.  Without the latter one cannot make general tensor products
of representations so that it is indeed remarkable that the formulation of 
geometric notions is possible.  This is what we outline, namely a gauge theory 
that has connections, principal bundles, associated bundles etc. using only
algebras and in particular not requiring anything from the theory of quantum groups. 

As such, the material in Section~2 should be rather more widely accessible 
than the coalgebra bundle version of the theory. In particular, 
it can be viewed as a critical first step towards a $C^*$ algebra or 
von Neumann algebra treatment. While  beyond our scope to actually consider 
operator theory and topological completions here (we work algebraically), it  
offers the possibility to link up with and extend other approaches to noncommutative
geometry based on $C^*$-algebras etc. We recall that in the $C^*$ algebra 
approach to noncommutative geometry, see \cite{Connes}, one traditionally 
works directly with vector bundles (as projective modules) and not principal 
bundles -- one would expect that the latter would
require some kind of group-like object such as a Hopf algebra but we see that 
this need not be the case. Also, although we do not develop a precise 
connection with the theory of subfactors at the present time, we note that 
our final data in terms of algebras is not unlike a subfactor inclusion. 
In that context one considers inclusions of von Neumann algebras with the larger one 
being viewed as some kind of `cross product' 
of the smaller one by some kind of  `paragroup'\cite{Ocn:str}. Similarly we show 
that if $A$ is an algebra acting on another 
algebra $P$  subject to a certain nondegeneracy condition then one can form a 
generalised `cross product' (which we call the `Galois product') of $P$ by $A$. 
In the subfactor case it is known that a special case corresponds to
some kind of (weak) Hopf algebra\cite{Boh:wea}, while similarly a special case 
of an algebra bundle corresponds to $A$ a Hopf algebra. The development of 
such an analogy on the one hand could provide a
gauge theory of subfactors (as well as a coalgebraic version of some aspects of their
theory) and on the other hand suggest the existence of a whole `Jones tower' 
of bundles. It would also connect with gauge theory from the point of view of algebraic
quantum field theory as in \cite{DopRob:why}\cite{FRS:sup} and many other works.

Linking up with $C^*$ or operator algebra results is the long-term motivation 
for the section. From a mathematical point of view, however, it should be 
stressed that our present results are strictly equivalent
to a subset of  the coalgebra bundle case. Part of the reformulation was
already hinted at in \cite{BrzMa:coa} where part of the data was expressed as an
algebra factorising into subalgebras $A,P$. The crucial exactness or
`Galois' condition in this form is what we provide now. It turns out to involve
traces over the vector space of $A$, which essentially forces us to 
finite-dimensional $A$. From this it is clear that the theory can be 
developed in two ways to cover the infinite-dimensional case: either
one needs to introduce operator completions which is the $C^*$ algebra 
or von Neumann algebra  direction mentioned above, or one needs to replace $A$ 
by its dual, a coalgebra, which then works for infinite-dimensional
coalgebras -- this is the approach taken in \cite{BrzMa:coa} and in the 
remaining sections of the present paper.

In Section~3 and Section~4 we continue with new results in the coalgebraic setting. We
provide the necessary formulation of associated bundles by
exploiting the recent work \cite{Brz:mod}. A small generalisation
of coalgebra bundles has been made in \cite{BrzHaj:coa} and we will
use in fact this formulation. Also, the notion of a connection which we
use here requires less structure than the one introduced in
\cite{BrzMa:coa}. In Section~5 we study frame
resolutions at this level. 

Finally, in Section~6, we show that the coalgebra theory allows one to include 
the crucial example of the monopole on the full 2-parameter family 
of Podle\'s quantum spheres\cite{Pod:sph}. Recall that
Podle\'s  classified all reasonable `quantum spheres' covariant
under the quantum group $SU_q(2)$, and until now the q-monopole has
been constructed\cite{BrzMa:gau} only for a diagonal
subfamily (the so-called standard quantum spheres). The general case 
requires the more general coalgebra bundle theory. The  bundle itself 
for all the quantum 2-spheres is in \cite{Brz:hom} and we provide on this 
now the required  connection. Similarly it is clear from their construction 
that all of the q-deformed symmetric spaces in the 
classification of \cite{NouSug:sym} should be constructable as coalgebra 
bundles, which includes the coalgebra bundle from which one would expect 
to project out a q-instanton. This is a second direction for further work.

We work algebraically over a general field $k$.  We use the usual
notations $\Delta c=c\o\tens c\t$ for a coproduct of a coalgebra
$C$ (summation understood). We also write $C^+=\ker\eps$ where $\eps$ 
is the counit. We write ${}_V\Delta(v)=v_{(1)}\tens v_{(\infty)}$ for a left
coaction on a vector space $V$, and $\Delta_V(v)=v_{(0)}\tens v\o$
for a right coaction. We also denote by $\Hom_A(V,W)$ the linear
maps commuting with a right action of an algebra $A$ and
by ${}_A\Hom(V,W)$ those commuting with a left action. Similarly,
$\Hom^C(V,W)$ for maps commuting with a right coaction of a
coalgebra $C$ and ${}^C\Hom(V,W)$ for a left coaction.  In general when 
we need to refer to the components of other elements $\chi^\#, \Psi(a\tens u)$ 
etc. of tensor product spaces we will use the upper bracket notation 
$\chi^\#=\chi^\#\uo\tens \chi^\#\ut$ etc. again with summation 
understood.

Also, we recall that for any algebra $P$, the universal 1-forms on
$P$ are $\Omega^1P=\ker\cdot_P\subset P\tens P$. The exterior
derivative $\extd: P\to \Omega^1P$ is $\extd u=1\tens u - u\tens 1$
for all $u\in P$. This extends to higher forms (see \cite{Kar:hom}
e.g.) $\Omega^{n}P \subset P^{\tens n+1}$ characterised by the
requirement that the products of all adjacent factors vanish, and $
\extd :\Omega^{n-1}P\to
\Omega^nP$,
\begin{equation}\label{extd}
\extd(u_1\tens\cdots \tens u_n)=\sum_{k=1}^{n+1} (-1)^{k-1} u_1\tens\cdots
\tens u_{k-1}\tens 1\tens
u_k\tens\cdots\tens u_n.
\end{equation}
With these definitions $\Omega P =
\bigoplus_{n=0}^\infty \Omega^nP$ is a graded differential algebra with
product given by juxtaposition and multiplication in $P$.

\section*{\normalsize\sc\centering 2. Galois actions and 
algebra factorisations}
\vspace{-.6\baselineskip}
\setcounter{section}{2}

Although we will continue mainly in an algebra-coalgebra setting in
later sections, we start with a more accessible version which
depends only on algebras and which should be useful for the
operator-algebraic version. We consider unital algebras and unital
algebra maps. An algebra factorisation means an algebra $X$ and
subalgebras $P,A$ such that the linear map $P\tens A\to X$ given by
the product in $X$ is an isomorphism.

\begin{prop}Cf. \cite{Tam:coe}\cite{Ma:book}\cite{CapSch:twi} algebra 
factorisations are in 1-1 correspondence with algebras
$P,A$ and $\Psi:A\tens P\to P\tens A$ such that
\[ \Psi(\cdot_A\tens\id)=(\id\tens\cdot_A)\Psi_{12}\Psi_{23},
\quad \Psi(1\tens u)=u\tens 1,
\quad\forall u\in P\]
\[ \Psi(\id\tens\cdot_P)=(\cdot_P\tens\id)\Psi_{23}\Psi_{12},
\quad \Psi(a\tens 1)=1\tens a,
\quad\forall a\in A.\]
In this case, given $e:A\to k$ a character, there is a left action
\[ \la:A\tens P\to P,\quad a\la u=(\id\tens e)\Psi(a\tens u),
\quad\forall a\in A,\ u\in P.\]
The subspace $M=P_e=\{m\in P|\ a\la m=e(a)m\ \forall a\in A\}$
forms a subalgebra.
\end{prop}
\proof Details of the stated equivalence are in \cite[pp. 299-300]{Ma:book}. Given 
$\Psi$ we define
$X=P\tens A$ with product
$(u\tens a)(v\tens b)=u\Psi(a\tens v)b$ for $u,v\in P$ and $a,b\in
A$. Given $X$ we define $\Psi$ by $au=\cdot_X\Psi(a\tens u)$. The
action $\la$ is also part of the proof in \cite{Ma:book} (where
$e=\eps$ the counit of a bialgebra). There is a similar right
action of $P$ on $A$ when $P$ is equipped with a character, which
we do not use.  From the point of view of $X$, $e$ on $A$ extends
to a left $P$-module map $e:X\to P$ obeying $e(au)=a\la u$ for all
$a\in A$ and $u\in P$.  Hence $M=\{u\in P|\ e(au)=ue(a)\ \forall
a\in A\}$, from which it is clear that $M$ is a subalgebra. One may
also see this from the equations for $\Psi$.
\endproof

Such factorisations are quite common. For example, they come up
naturally as part of Hopf algebra factorisations\cite{Ma:phy}\cite{Ma:book}. 
Another example
is the braided tensor product $A\und\tens B$ of two algebras, see
\cite{Ma:book}. In our geometrical picture, $P$ plays the role of
the algebra of functions of the `total space' of a principal
bundle, and $A$ plays the role of the group algebra of the
structure group. The subalgebra $M$ plays the role of the functions
on the base. The algebra $X$ is not usually considered but plays
the role of the cross product $C^*$-algebra of the functions on the
total space by the action of the structure group.

\begin{prop} In the setting above, the map $\tilde\chi:A\tens P\tens P\to P$
defined by $\tilde\chi(a\tens u\tens v)=(a\la u)v$ descends to
$\chi:A\tens P\tens_M P\to P$. We say that the factorisation is
{\em Galois} if there exists $\chi^\#:P\to P\tens_M P\tens A$ such
that
\[ \tr_A(\chi^\#\circ\chi)=\id_{P\tens_M P},\quad
(\chi\tens\id_A)(\id_A\tens\chi^\#)=\tau:A\tens P\to P\tens A\]
where $\tau$ is the usual flip or transposition map. We call $P(M,A,\Psi,e)$
a {\em copointed algebra bundle}.
\end{prop}
\proof We have $a\la(um)=e(aum)=e(u_ia^im)=u_ie(a^im)=u_im e(a^i)=(a\la u)m$
for all $a\in A$, $u\in P$ and $m\in M$, as required. Here
$\Psi(a\tens u)=u_i\tens a^i$ is a notation (sum over $i$). The
rest is a definition. This can also be obtained from the $\Psi$
equations. \endproof

This is the analogue of the Galois condition in \cite{BrzMa:coa},
which in turn is motivated from the theory of quantum principal
bundles and, independently, the theory of Hopf-Galois extensions in
the Hopf algebra case. In geometrical terms the map $\chi$ plays
the role of the action of the Lie algebra $\cg$ of the structure
group of a principal bundle on its algebra of functions: if
$\xi\in\cg$ one has a left-invariant vector field $\tilde\xi$ given
by differentiating the action corresponding to $\la$. The element
$\chi^\#=\chi^\#(1)\in P\tens_MP\tens A$ is particularly important and
plays the role of the `translation map' of the principal bundle.

Notice, however, that a factorisation can be Galois
only if $A$ is finite-dimensional. This should not unduly worry us 
since our formulation is mainly intended as a precursor to an 
operator-theoretic treatment where infinite-dimensional $A$ would 
be allowed subject to topological completions and
trace class conditions. To avoid all that in the infinite-dimensional 
case one should of course use the coalgebra formulation as in 
later sections. Meanwhile, let us note that even finite-dimensional 
$A$ are not uninteresting -- the algebra $P$ and the factorising
algebra can in principle both be infinite-dimensional. A similar 
situation pertains with the theory of subfactors where the two von 
Neumann algebras are typically infinite-dimensional
but the case where their `ratio' is in some sense finite is 
still very interesting.

There is an obvious notion of a $\Psi$-module associated to an
algebra factorisation, namely a left $P$ module and $A$ module $V$ such
that
\eqn{Psimod}{ a\la(u\la v)=\la\circ(\Psi(a\tens u)\la v),\quad
\forall a\in A,\ u\in P, \ v\in V.}
Explicitly we require $a\la(u\la v) = u_i\la(a^i\la v)$,where
$\Psi(a\otimes u) = u_i\otimes a^i$.
This is what corresponds to a left $X$-module. This point of view
suggests a natural slight generalisation  of the above, replacing
$e$ by the requirement of a map $\tilde e:A\to P$.

\begin{prop} Let $(P,A,\Psi)$ be an algebra factorisation datum
as in Proposition~2.1.
Linear maps $\tilde e:A\to P$ such that
\[ \tilde e(ab)=\Psi(a\tens\tilde e(b))\uo\tilde e(\Psi(a\tens
\tilde e(b))\ut),\quad \tilde e(1)=1,\quad\forall a,b\in A,\]
are in 1-1 correspondence with extensions of the left regular
action of $P$ on itself to a $\Psi$-module structure on $P$.
Given $\tilde e$, we define
\[ a\la u=\Psi(a\tens u)\uo\tilde e(\Psi(a\tens u)\ut),\quad
\forall a\in A,\ u\in P\]
and conversely, given such an extension, we set $\tilde e(a)=a\la
1$. In this situation the space
\[ M=M_{\tilde e}=\{m\in P| \ a\la m=\tilde e(a)m,\ \forall a\in A\},\]
is a subalgebra of $P$ and $\tilde\chi$ as in Proposition~2.2
descends to a map $\chi$.\end{prop}
\proof We define the linear map $\la:A\tens P\to P$ as stated
and verify first equation (\ref{Psimod}) as
\[\Psi(a\tens uv)\uo\tilde e(\Psi(a\tens uv)\ut)
=u_i\Psi(a^i\tens v)\uo\tilde e(\Psi(a^i\tens v)\ut)=\Psi(a\tens u)
\uo(\Psi(a\tens u)\ut\la v),\]
where we used the second of factorisation properties in
Proposition~2.1. We also used the shorthand $\Psi(a\tens u)=u_i\tens a^i$
as before. Next, we check that $\la$ is indeed an action,
\align{(ab)\la v\equad &&=\Psi(ab\tens u)\uo\tilde e(\Psi(ab\tens u)\ut)=
\Psi(a\tens u_i)\uo\tilde e(\Psi(a\tens u_i)\ut b^i)\\
&&=\Psi(a\tens u_i)\uo (\Psi(a\tens u_i)\ut\la\tilde e(b^i))
=a\la(u_i\tilde e(b^i))=a\la(b\la u)}
as required. We used the first of the factorisation properties of
$\Psi$ and the assumed condition on $\tilde e$, which can be written
as $\tilde e(ab)=a\la\tilde e(b)$ in terms of $\la$. We then used
(\ref{Psimod}) already proven. Finally, $1\la u=u\tilde e(1)=1$ so
$\la$ is indeed an action. Conversely, given an action $\la$ making
$P$ a $\Psi$-module we define $\tilde e(a)=a\la 1$. Then $\tilde e(ab)=
a\la(b\la 1)=a\la \tilde e(b)$ and $a\la u=a\la (u\la 1)=u_i (a^i\la 1)
=u_i\tilde e(a^i)$ (using (\ref{Psimod})), as required. The remaining
facts follow easily from the definition of $M$. It can also be characterised
equivalently as
\eqn{newM}{ M=\{m\in P|\ a\la(um)=(a\la u)m,\ \forall u\in P,\ a\in A\}}
in view of (\ref{Psimod}) and the definition of $\la$. \endproof

We note that

\begin{lem} In the setting of Proposition~2.3, for any $\Psi$-module $V$
there is a natural notion of `invariant' subspace
\eqn{Vinv}{V_0=\{v\in V|\ a\la v=\tilde e(a)\la v,\quad \forall 
a\in A\}\subseteq V}
which is a left $M$-module by restriction of the action of $P$.
\end{lem}
\proof For all $a\in A, m\in M$ and $v\in V_0$, we have
$a\la(m\la v)=\la\circ(\Psi(a\tens m)\la v)=m_i\la(\tilde e(a^i)\la
)=(m_i\tilde e(a^i))\la v=(a\la m)\la v=(\tilde e(a)m)\la v=\tilde
e(a)\la(m\la v)$, so $m\la v\in V_0$ as well. \endproof

The subalgebra $M$ itself is a case of such an invariant subspace.
When there is a corresponding $\chi^\#$, we call $P(M,A,\Psi,\tilde
e)$ an algebra bundle. The copointed case is $\tilde e(a)=e(a)1$.
The construction has a natural converse.

\begin{lem} In an algebra bundle, $\chi^\#=\chi^\#(1) $ obeys

(a) $\chi^\#\, a=a\la\chi^\#$ for all $a\in A$, where the action on
$P\tens_MP\otimes A$ is on its first factor.

(b) $\chi^\#(uv)=\chi^\#(u)\uo\, v\tens \chi^\#(u)\ut$ for all $u,v\in P$

(c) $\chi^\#\uo\,  (\chi^\#\ut\la u)=u\tens_M1$.

\noindent Here $\chi^\#\uo\in P\tens_MP$ and $\chi^\#\ut\in A$ 
for all $u\in P$.
\end{lem}
\proof From its definition, it is evident that
\[ \chi(ab\tens u\tens v)=(ab\la u)v=\chi(a\tens b\la u\tens v),\quad
\chi(a\tens u\tens vw)=\chi(a\tens u\tens v)w\]
for all $u,v,w\in P$ and $a,b\in A$. Parts (a) and (b) are just the
corresponding properties for
$\chi^\#$. Thus,
\align{a\la\chi^\#\equad&&=(\trace_A\chi^\#\circ\chi\tens\id)(a\la\chi^\#)
=(\id\tens f^a)\chi^\#(\chi(e_a\tens a\la \chi^\#\uo))\tens\chi^\#\ut\\
&&=(\id\tens f^a)\chi^\#\circ\chi(e_a a\tens \chi^\#\uo)\tens \chi^\#\ut
=(\id\tens f^a)\chi^\#(1)\tens e_aa=\chi^\#a}
where $\{e_a\}$ is a basis of $A$ and $\{f^a\}$ a dual basis. Similarly,
\align{ \chi^\#(u)\uo v\tens \chi^\#(u)\ut\equad &&
=(\trace_A\chi^\#\circ\chi(\chi^\#(u)\uo v))
\tens \chi^\#(u)\ut\\
&&=(\id\tens f^a)\chi^\#(\chi(e^a\tens \chi^\#(u)\uo)v)\tens \chi^\#(u)\ut
=\chi^\#(uv).}
We then deduce part (c) from part (b) as
\align{ \chi^\#(1)\uo (\chi^\#(1)\ut\la u)\equad&&
=\chi^\#(1)\uo(e_a\la u)\<f^a,\chi^\#(1)\ut\>
\\
&&=(\id\tens f^a)\chi^\#(1.(e_a\la u))
=(\id\tens f^a)\chi^\#\circ\chi(e_a\tens u)=u\tens 1.}
\endproof

These correspond to important properties of the translation map in 
differential
geometry derived in the Hopf algebraic setting in \cite{Brz:tra}\cite{Sch:rep}.

\begin{theorem} Let $P,A$ be algebras and $P$ a left $A$-module under
an action $\la$.
We define $M$ by (\ref{newM}) and say that the action is Galois if 
 $\chi$ defined as in Proposition~2.2 has
a corresponding $\chi^\#$. In this case there exists a unique algebra
factorisation $X= P\und\otimes_\Psi A$ such that $P,A$ form an algebra 
bundle and $P$ is a $\Psi$-module (cf.\ eq.~(\ref{Psimod})) via product
in $P$ and the action $\la$. Explicitly, 
\[ \Psi(a\tens u)=\chi(a\tens u\chi^\#\uo)\tens \chi^\#\ut,
\quad \tilde e(a)=a\la 1,\quad \forall a\in A,\ u\in P.\]
We call the corresponding algebra factorisation $X=P\und\tens_\Psi A$ the
{\em Galois product} associated to a Galois action of an algebra $A$
on an algebra $P$.
\end{theorem}
\proof Here we define $M$ and $\tilde\chi$ directly from the action $\la$; it is easy
to see that $M$ is a subalgebra and that $\tilde\chi$ descends to a map $\chi$. 
We assume the existence of a corresponding $\chi^\#$ obeying the conditions
in Proposition~2.2. For the purposes of this proof, we now write
 $\chi^\#=\chi^\#\uo\tens_M\chi^\#\ut\tens\chi^\#\uth$ (a more explicit
notation than the one before) and we let $\chi'^\#$ be a second
copy of $\chi^\#$. Then the map $\Psi$ explicitly reads
$$
\Psi(a\otimes u) = (a\la(u\chi^\#\uo))\chi^\#\ut\otimes \chi^\#\uth
$$
and we have,
\align{\equad&& (\id\tens\cdot_A)\Psi_{12}\Psi_{23}(a\tens b\tens u)=
(\id\tens\cdot_A)\Psi_{12}(a\tens (b\la(u\chi^\#\uo))\chi^\#\ut)
\tens\chi^\#\uth\\
&&=\left(a\la\left((b\la(u\chi^\#\uo))\chi^\#\ut \chi'^\#\uo\right)\right)
\chi'^\#\ut\tens
\chi'^\#\uth \chi^\#\uth\\
&&= \left(a\la\left((b\la(u\chi^\#\uo))\chi^\#\ut (\chi^\#\uth\la
\chi'^\#\uo)\right)\right)\chi'^\#\ut\tens
\chi'^\#\uth\\
&&= (a\la(b\la(u\chi'^\#\uo)))\chi'^\#\ut \chi'^\#\uth=\Psi(ab\tens u)}
using parts (a) and then (c) of the lemma and that $\la$ is an action.
On the other side, we have
\align{(\cdot_P\tens\id)\equad && \Psi_{23}\Psi_{12}(a\tens u\tens v)
=(\cdot_P\tens\id)\Psi_{23}((a\la(u\chi^\#\uo))\chi^\#\ut\tens
\chi^\#\uth\tens v)\\
&&=(a\la(u\chi^\#\uo))\chi^\#\ut (\chi^\#\uth\la(v\chi'^\#\uo))
\chi'^\#\ut\tens\chi'^\#\uth\\
&&=(a\la(u v \chi'^\#\uo))\chi'^\#\ut\tens\chi'^\#\uth=\Psi(a\tens uv)}
using part (c) of the lemma. The computations for $\Psi(a\tens 1)$
and $\Psi(1\tens u)$ are more trivial and left to the reader. We need
\[ \chi^\#\uo\chi^\#\ut\tens\chi^\#\uth
=\chi(1\tens\chi^\#\uo\tens \chi^\#\ut)\tens\chi^\#\uth=1\tens 1\]
for the latter case.

Hence we have a factorisation datum and by Proposition~2.1 we
have an algebra $X$ built on
$P\tens A$ with the cross relations $(1\tens a)(u\tens 1)=\Psi(a\tens u)$.
We now define $\tilde e(a)=a\la 1$ and check easily that
$\tilde e(ab)=a\la\tilde e(b)$ as required, and that
$\Psi(a\tens u)\uo\tilde e(\Psi(a\tens u)\ut)=(a\la (u\la \chi^\#\uo))
\chi^\#\ut (\chi^\#\uth\la 1)
=a\la u$ by part (c) of the lemma. Hence is the converse to the
preceding proposition.

To prove that $P$ is a $\Psi$-module, we take any $a\in A$, $u,v\in P$
and use the explicit form of $\Psi$ above and part (c) of the lemma 
to compute
$$
\cdot\circ(\Psi(a\tens u)\la v) =
(a\la(u\chi^\#\uo))\chi^\#\ut(\chi^\#\uth\la v) = a\la(uv).
$$
Finally, suppose there is another factorisation $\Psi'$ such that $P$ is
a $\Psi'$-module, and let $\Psi'(a\otimes u)= u_i\otimes a^i$ for all
$a\in A$, $u\in P$. Then
$$
 \Psi(a\otimes u) = (a\la(u\chi^\#\uo))\chi^\#\ut\otimes \chi^\#\uth =
u_i(a^i\la \chi^\#\uo)\chi^\#\ut\otimes \chi^\#\uth =u_i\otimes a^i = 
\Psi'(a\otimes u),
$$
where we used the definition of $\chi^\#$. This proves the uniqueness of
$\Psi$.  
\endproof

\begin{example} Let $q$ be a primitive $n$'th root of 1. The $n\times n$ matrices
factorise as $M_n(\C)=\C\Z_n\cdot\C\Z_n$, where the two copies of
$\Z_n$ are generated by \[ g=\pmatrix{1&0&0\cdots 0\cr0&q&0\cdots
0\cr&\vdots&\vdots&\cr0&\cdots&0&q^{n-1}},\quad
h=\pmatrix{0&1&0&0\cdots0\cr0&0&1&0\cdots
0\cr&\vdots&&\vdots&\cr1&0&0\cdots&0}\] obeying $hg=qgh$, so
\[\Psi(h\tens g)=qg\tens h,\quad  \Psi(1\tens g)=g\tens 1,\quad
 \Psi(h\tens 1)=1\tens h,\quad\Psi(1\tens 1)=1\tens 1.\]
 The nontrivial character $e(h)=q$ gives
 \[ h\la g^m=q^{m+1}g^m\]
 and hence $M=\C 1$. The result is Galois, with
\[ \chi(h^m\tens g^k\tens g^l)=q^{m(k+1)}g^{k+l},\quad\chi^\#(g^m)=n^{-1}
\sum_{a,b}q^{-ab}g^{b-1}\tens g^{m-b+1}\tens h^a.\]
 \end{example}
\proof We identify $A=\C\Z_n=\C[h]/h^n=1$ and $P=\C\Z_n=\C[g]/g^n=1$ as the
two subalgebras. The relations $hg=(1\tens h)(g\tens 1)=
\Psi(h\tens g)=q(g\tens 1)(1\tens h)$ give the form of $\Psi$.
This extends uniquely to a solution of the factorisation equations
in Proposition~2.1 as $\Psi(h^m\tens g^k)=q^{mk}g^k\tens h^m$. Actually,
this is an example of a braided tensor product $M_n(\C)=\C\Z_n\und\tens\C\Z_n$
in the braided category of anyonic or $\Z_n$-graded spaces. The character $e$
then gives the action shown as $h\la g^m=q^m g^m e(h)$. Hence
$\sum_m a_m g^m\in M$ iff $a_m(q^{m+1}-q)=0$ for all $m$, which
means $M=\C1$. We also obtain $\chi$ as shown and one may verify
that $\chi^\#$ as stated fulfills the requirements in
Proposition~2.2. \endproof

In this example $A$ is actually a Hopf algebra and 
$\tilde e(h)=q1$ as here yields a bundle with is equivalent (in the
coalgebra bundle version) to a Hopf algebra bundle as in
\cite{BrzMa:gau}. On the other hand, other choices of $\tilde e$
yield algebra bundles which are not equivalent to Hopf algebra
bundles, i.e. strict examples of our more general theory. We
examine the $n=2$ case in detail:

\begin{example} The factorisation $M_2(\C)=\C\Z_2\cdot\C\Z_2$ as above
(with $q=-1$) admits a family of algebra bundle structures
parametrized by $\theta\in [0,2\pi)$, with
\[ \tilde e(h)=\cos(\theta)+\imath g\sin(\theta).\]
The associated Galois action is
\[ h\la g^k=(-1)^k g^k (\cos(\theta)+\imath g\sin(\theta)).\]
\end{example}
\proof We have $A=\C[h]/h^2=1$ and $P=\C[g]/g^2=1$. We require $\tilde e$
of the form
\[ \tilde e(1)=1,\quad \tilde e(h)=\alpha+\imath\beta g\]
(say) obeying the condition in Proposition~2.3. The non-empty case
is $1=\tilde e(1)=\tilde e(h.h)=\Psi(h\tens\tilde e(h))\uo\tilde
e[\Psi(h\tens\tilde e(h))\ut]=\alpha\tilde e(h)-\imath \beta
g\tilde e(h)=(\alpha-\imath \beta g)(\alpha+\imath \beta
g)=\alpha^2+\beta^2$. This admits many solutions over $\C$, a
natural family being those where $\alpha,\beta$ are real, i.e. on a
circle parametrized by $\theta$. On the other hand, in the case
equivalent to a Hopf algebra bundle, $P$ would be an $A$-module
algebra. This happens when $(h\la g)^2=\cos^2(\theta)
-\sin^2(\theta)+2\imath\sin(\theta)\cos(\theta)g=1$, which is
 exactly when $\theta=0,\pi$. The first case is trivial and the second is
 the $n=2$ case of the preceding Example~2.6.

Next, we consider $m=a+bg$ such that $h\la m=\tilde e(h)m$. It is easy to see 
that this happens iff $b=0$,
provided $\sin(\theta)\ne 0$ or $\cos(\theta)\ne 0$ (one of which
is always the case). Hence $M=\C 1$. Finally, we have
$\chi(h^m\tens g^k\tens g^l)=g^{k+l}((-1)^k(\cos(\theta)+\imath
\sin(\theta)g))^m$ which we can write as a map $P\tens P\to
A^*\tens P$. Identifying $A^*=\C\Z_2$ with generator $c$, say, the
map is
\[ g^k\tens g^l\mapsto c_++c_-(-1)^k\cos(\theta)\tens g^{k+l}
+c_-\imath(-1)^k\sin(\theta)\tens g^{k+l+1}\] where $c_\pm=(1\pm
c)/2$. (This is the map $\chi$ in the corresponding coalgebra
bundle version). Invertibility of this map is equivalent to the
existence of $\chi^\#$ in the present setting; in fact the map has
determinant 1 in the obvious basis $\{g^k\tens g^l\}$ and
$\{c^k\tens g^l\}$ and is therefore invertible. \endproof

The corresponding factorisation over $\R$ is the quaternion algebra
and provides a counterexample to the existence of $\tilde e$:

\begin{example} Over $\R$, the quaternion algebra $\H=\span\{1,i,j,k\}$
obeying $i^2=j^2=k^2=-1$
and $ij=k$ etc., is a factorisation $\H=\R[i]\R[j]$ where
$\R[i]=\C$ as a 2-dimensional algebra over $\R$. One has
\[ \Psi(j\tens i)=-i\tens j,\quad \Psi(1\tens i)=i\tens 1,\quad
 \Psi(j\tens 1)=1\tens j,\quad\Psi(1\tens 1)=1\tens 1.\]
This factorisation admits no map $\tilde e$.
\end{example}
\proof The factorisation is evident, with $P=\R[i]$  and $A=\R[j]$
(the quotient of
polynomials by the relation $i^2=-1$ and $j^2=-1$ respectively). Now suppose
a linear map $\tilde e:\R[j]\to\R[i]$ of the form
\[ \tilde e(1)=1, \quad \tilde e(j)=\alpha+i\beta.\]
Then a similar computation to the one above yields this time
$-1=\tilde e(-1)=\tilde e(j.j)
=\alpha^2+\beta^2$, which has no solutions over $\R$. \endproof

Returning to the general theory,

\begin{prop} An algebra bundle is {\em trivial} or `cleft' if there is an
invertible element $\Phi=\Phi\uo\tens\Phi\ut\in P\tens A^{\rm op}$ such that
\[ \Phi a=a\la \Phi,\quad\forall a\in A,\]
where the product from the right is in $A$. In this case, $P\isom
\Hom_k(A,M)$ as a left $A$-module and right $M$-module.

Moreover, $(P,A,\Psi,\tilde e)$ in Proposition~2.3 is a trivial (cleft)
algebra bundle if there exists an invertible $\Phi\in P\tens A^{\rm op}$ 
obeying the above condition,  with
\[ \chi^\#(u)=\Phi\uo\tens_M\Phi\umo u\tens\Phi\umt\Phi\ut,\quad 
\forall u\in P.\]
\end{prop}
\proof The isom $\Theta:\Hom_k(A,M)\to P$ is
\[ \Theta(f)=\Phi\uo f(\Phi\ut),\quad \Theta^{-1}(u)(a)=\Phi\umo
((\Phi\umt a)\la u).\]
Here $\Theta$ is a left $A$-module map since the image of $f$ is 
in $M$ and
$\Phi$ obeys the condition above. Next, the latter condition is 
equivalent to the condition
\eqn{Phiinv}{ \Psi(a\tens\Phi\umo)\Phi\umt=\tilde e(a)\Phi^{-1}}
for $\Phi^{-1}$ (just compute $\tilde e(a)\tens 1=a\la 1\tens 1
=a\la(\Phi\umo\Phi\uo)\tens\Phi\ut \Phi\umt$ using (\ref{Psimod})). Hence
$a\la(\Phi\umo((\Phi\umt\la u)))=\cdot\circ\Psi(a\tens\Phi\umo)\Phi\umt\la u
=\tilde e(a)\Phi\umo((\Phi\umt\la u)$, i.e. $\Phi\umo((\Phi\umt\la u\in M$ for
all $u\in P$. In particular, this implies that $\Theta^{-1}(u):A\to M$ as
required. This then provides the required inverse since
$\Theta\circ\Theta^{-1}(u)=\Phi\uo(\Theta^{-1}(u))(\Phi\ut)=
\Phi\uo\Phi\umo((\Phi\umt\Phi\ut)\la u)=u$ from the definitions, and
$(\Theta^{-1}\circ\Theta(f))(a)=\Phi\umo(\Phi\umt a\la\Theta(f))
=\Phi\umo\Theta(\Phi\umt a\la f)=\Phi\umo\Phi\uo f(\Phi\ut\Phi\umt a)=f(a)$
by the left $A$-module property of $\Theta$.

For the second part, given a factorisation datum and $\tilde e$, 
we define $\chi^\#$ as shown. Then
\[ \chi(a\tens\chi^\#(u)\uo)\tens\chi^\#(u)\ut=(a\la\Phi\uo)\Phi\umt u
\tens\Phi\umt\Phi\ut=\Phi\uo\Phi\umo u\tens\Phi\umt\Phi\ut a=u\tens a\]
by the property of $\Phi$. On the other side,
\[\trace_A\chi^\#\circ\chi(u\tens_M v)=\Phi\uo\tens_M\Phi\umo
(\Phi\umt\Phi\ut\la u)v=\Phi\uo\Phi\umo(\Phi\umt\Phi\ut\la u)\tens_Mv=u\tens_Mv\]
since $\Phi\umo(\Phi\umt \la (\Phi\ut\la u))\in M$ by the same proof as above.
Hence the action of $A$ on $P$ is Galois in this case.
\endproof 

\begin{prop} A bundle automorphism is an invertible linear map
$P\to P$ which is a right
$M$-module map and a left $A$-module map. The group of bundle
automorphisms is in correspondence with invertible $f\in P\tens
A^{\rm op}$ such that
\[ \Psi(a\tens f\uo) f\ut=f\uo\tens f\ut a,\quad \forall a\in A.\]
When the bundle is trivial, such $f$ correspond to 
invertible elements $\gamma\in M\tens A^{\rm op}$ by $f=\Phi\gamma\Phi^{-1}$
multiplied in $P\tens A^{\rm op}$.
\end{prop}
\proof Note first of all that the set of such elements in
$P\tens A^{\rm op}$ form
a group. Thus, $\Psi(a\tens f\uo g\uo)g\ut
f\ut=f\uo{}_i\Psi(a^i\tens g\uo)g\ut f\ut=f\uo{}_i g\uo\tens g\uo
a^i f\ut=f\uo g\uo\tens g\ut f\ut a$ when $f,g$ obey this
condition. The relation between such $f$ and automorphisms $F$ is
\[F(u)=f\uo(f\ut\la u),\quad f=F(\chi^\#\uo)\chi^\#\ut\tens \chi^\#\uth.\]
Thus, given $f$ it is evident that $a\la F(u)=a\la( f\uo(f\ut\la u))
=\cdot\circ\Psi(a\tens f\uo)f\ut\la u=f\uo(f\ut a\la u)=F(a\la u)$ by
(\ref{Psimod}) and the property of $f$, so $F$ is a left $A$-module map 
(it is clearly a right $M$-module map as well). Also from this, it is
immediate that the product in $P\tens A^{\rm op}$ maps over to the
composition of bundle transformations. Finally, the inverse of the 
construction is as shown using the properties of $\chi^\#$. Thus,
$F(\chi^\#\uo)\chi^\#\ut(\chi^\#\uth\la u)=F(u)$ by Lemma~2.5(c), and
when $F$ is defined by $f$, the inversion formula yields
$(f\uo\la(f\ut\la\chi^\#\uo))\chi^\#\ut\tens \chi^\#\uth=f\uo\chi(f\ut\tens
\chi^\#\uo\tens\chi^\#\ut)\tens\chi^\#\uth=f$ from the definition of 
$\chi^\#$. 

In the case of a trivial bundle, we define $f$ as shown and verify
\align{\Psi(a\tens f\uo)f\ut\equad &&=\Psi(a\tens\Phi\uo\gamma\uo
\Phi\umo)\Phi\umt\gamma\ut\Phi\ut\equad \\
&&=(\Phi\uo\gamma\uo){}_i\Psi(a^i\tens\Phi\umo)\Phi\umt\gamma\ut\Phi\ut\\
&&=(a\la(\Phi\uo\gamma\uo))\Phi\umo\tens \Phi\umt\gamma\ut\Phi\ut
=f\uo\tens f\ut a}
using the properties of $\Psi$, (\ref{Phiinv}) and that $\gamma\uo\in M$. 
Conversely, given $f$ we define $\gamma=\Phi^{-1}f\Phi$ (product in 
$P\tens A^{\rm op}$ and verify using (\ref{Phiinv}) that $a\la\gamma
=\tilde e(a)\gamma$ so that $\gamma\in M\tens A^{\rm op}$. \endproof

Next, even though $A$ is only an algebra, its action on $P$ extends
naturally to tensor powers and hence to the universal exterior
differentials $\Omega^n P\subset P^{\tens (n+1)}$.

\begin{prop} In the setting of Proposition~2.3, $\Omega^nP$ is a
$\Psi^\bullet$-module with
$P$ acting by left multiplication and
\[ a\la (u_0\tens\cdots\tens u_n)=\Psi^{n+1}(a\tens u_0\tens
\cdots\tens u_n)
\uo\tilde e[\Psi^{n+1}(a\tens u_0\tens\cdots\tens u_n)\ut]\]
where
$\Psi^\bullet|_{\Omega^nP}=\Psi^{n+1}=\Psi_{n,n+1}\cdots\Psi_{12}$
defines another factorisation datum $(\Omega P,\Psi^\bullet, A)$.
It is such that $\Psi^\bullet\circ(\id\tens\extd)
=(\extd\tens\id)\circ\Psi^\bullet$.
\end{prop}
\proof That $(\Omega P,\Psi^\bullet,A)$ is another factorisation datum
is an elementary proof by induction repeatedly using the
factorisation properties in Proposition~2.1 and the product in
$\Omega P$ (which is just inherited from the product in $P$); it is
left to the reader. Applying Proposition~2.3 to this new
factorisation, with $\tilde e:A\to P\subset\Omega P$ then gives a
$\Psi^\bullet$-module. We then restrict the action to ones of
$A,P$. \endproof

Armed with this, we can define a connection as an equivariant splitting
of $\Omega^1P\supseteq P(\Omega^1M)P$ as in
\cite{BrzMa:gau}\cite{BrzMa:coa}. More precisely, we require that
$\Pi$ has kernel $P(\Omega^1M)P$, is a right $P$-module map and $\Pi\circ\extd$ 
is left $A$-module map. 
Such projections turn out to be in
1-1 correspondence with $\omega\in \Omega^1P\tens A$ such that

i) $\omega\uo\tilde e(\omega\ut)=0$

ii) $\tilde\chi(a\tens\omega\uo)\tens\omega\ut=1\tens a-\tilde e(a)\tens 1$

iii) $\omega a=\Psi^\bullet(a\tens\omega\uo)\omega\ut$

Here the correspondence is
\[ \Pi(u\tens v)=\omega\uo\tilde\chi(\omega\ut\tens u\tens v)\]
(using $\chi^\#$ in the reverse direction). We will provide this in more 
detail in the next
section in the coalgebra setting.

There is also a theory of associated bundles. In fact, one has and needs 
two kinds of
associated bundles; given an algebra bundle and a right
$A$-module $V_R$ we have
\[ E=\{\sum_k v_k\tens u_k\in V_R\tens P|\ \sum_k v_k\ra a\tens u_k
=\sum_k v_k\tens a\la u_k,\ \forall
a\in A\}\subseteq V_R\tens P\]
as a natural right $M$-module by right multiplication in $P$. And given
a left $A$-module $V_L$ we have
\[ \bar E=\{\sum_k u_k\tens v_k\in P\tens V_L|\ 
\sum_k\Psi(a\tens u_k)\la v_k=\sum_k \tilde e(a)u_k\tens 
v_k,\ \forall a\in A\}
\subseteq P\tens V_L\]
as a natural left $M$-module. Here the $\bar E$ is the natural
`invariant' subspace from Lemma~2.4 for the $\Psi$-module structure
of $P\tens V_L$ provided by the following lemma.

\begin{lem} If $V$ is a left $A$-module then $P\tens V$ is a
$\Psi$-module where $P$ acts by multiplication from the left and $A$ acts by
\[ a\la (u\tens v)=\Psi(a\tens u)\la v.\]
\end{lem}
\proof We check first that $A$ acts as shown. Thus,
$(ab)\la(u\tens v)=\Psi(ab\tens u)\la v=
\Psi(a\tens u_i)\la (b^i\la v)=\Psi(a\tens (b\la(u\tens v))\uo)\la
(b\la(u\tens v))\ut
=a\la(b\la(u\tens v))$ using the definitions and the factorisation
 property of $\Psi$.
Here $b\la(u\tens v)=(b\la (u\tens v))\uo\tens (b\la(u\tens v))\ut$
is a notation. This then forms a $\Psi$-module since
\[ a\la(uu'\tens v)=\Psi(a\tens uu')\la v=u_i\Psi(a^i\tens u')\la v
=u_i(a^i\la(u'\tens v))=\la\circ(\Psi(a\tens u)\la(u'\tens v))\]
as required.
\endproof

Sections of these bundles are $M$-valued $M$-module maps from $E$, 
$\bar E$ respectively. When $P$ is flat 
over $M$ and $\Psi$ has a certain adjoint $\Psi^\#$, one can show that 
\[ \Hom_A(V_L,P)\isom {}_M\Hom(\bar E,M),\quad \Hom(V_R,P)_0\isom\Hom_M(E,M)\]
as right $M$-modules, left $M$-modules respectively. In the first case, 
if $\phi\in\Hom_A(V_L,P)$ then the corresponding section of $\bar E$ is
${\bar s}_{\phi}(u\tens v)=u\phi(v)$. In the second case, $\Hom(V_R,P)$ is a 
left $\Psi$-module in a similar manner to Lemma~2.13
(coinciding with it in the finite-dimensional case, namely $(a\la \phi)(v)
= \trace_A\Psi(a\tens \phi(v\la(\ )))$. If $\phi\in\Hom(V_R,P)_0$ then 
the corresponding
section of $E$ is $s_{\phi}(v\tens u)
=\phi(v)u$.  The proof of these assertions will be given in Section~4 in 
the coalgebra setting with $\chi^{-1}$ and $\psi^{-1}$ in the roles of 
$\chi^\#$ and $\Psi^\#$.

When $V_L$ and $V_R$ are finite-dimensional then
\[ E=\Hom_A(V_R^*,P),\quad \bar E=\Hom(V_L^*,P)_0,\]
so that each bundle can be viewed as the space of sections of the
other. Moreover, the constructions generalise
directly to form-valued sections by using $\Psi^\bullet$ in place of $\Psi$.
One may then proceed to frame bundles etc. Thus, one has a
covariant derivative
\[ \nabla: E\to E\tens_M\Omega^1M,\quad \bar\nabla:\bar E\to
\Omega^1M\tens_M\bar E\]
associated to a suitable (strong) connection in the pointed case. 
By definition a frame resolution
is an associated bundle equipped with a canonical form such that 
$E\isom\Omega^1M$, and in
this case $\nabla$ plays the role of Levi-Civita connection etc, 
along the lines in \cite{Ma:rie}. This and the rest of the theory will be
provided in Section~4, in our preferred coalgebra bundle setting.

Finally, we give the situation in the case of trivial (cleft) algebra bundles.
In this case sections correspond to `matter fields' on the base $M$,
\[ \Hom_A(V_L,P)\isom \Hom_k(V_L,M),\quad \Hom(V_R,P)_0=\Hom_k(V_R,M).\]
The first isomorphism sends $\bar f\in\Hom_k(V_L,M)$ to the map 
${\phi}_{\bar f}(v)=\Phi\uo\bar f(\Phi\ut\la v)$. The second isomorphism sends
$f\in\Hom_k(V_R,M)$ to ${\phi}_f(v)=f(v\ra\Phi\umt)\Phi\umo$. Similarly,
(strong) connections $\omega$ are determined by `gauge fields' 
$\alpha\in\Omega^1M\tens A$ such that $\alpha\uo\tilde e(\alpha\ut)=0$, 
according to
\[ \omega=\Phi\uo\alpha\uo\Phi\umo\tens \Phi\umt\alpha\ut\Phi\ut+
\Phi\uo\extd \Phi\umo\tens\Phi\umt\Phi\ut.\]
Proofs will again be given in the following sections, in the coalgebra
setting. The covariant derivative on these matter fields and their 
gauge transformation by 
$\gamma\in M\tens A^{\rm op}$ then take on the familiar form for algebraic 
gauge theory on trivial bundles (see \cite[Sec. 3]{Ma:rem}).

\section*{\normalsize\sc\centering 3. Coalgebra bundles and 
connections}\vspace{-.6\baselineskip}
\setcounter{section}{3}
\setcounter{definition}{0}

We switch now to the coalgebra version of the theory, where $A$ is
replaced by a coalgebra $C$. This is the original theory of
coalgebra bundles\cite{BrzMa:coa}, which we extend further. The
coalgebra version involves less familiar notations but has
advantages in a purely algebraic treatment.

\begin{definition}\label{ent}\cite{BrzMa:coa}
A coalgebra $C$ and algebra $P$ are {\em entwined} by $\psi: C\tens
P\to P\tens C$ if
\begin{equation}\label{diag.A}
\psi\circ(\id\tens \cdot) = (\cdot\tens \id)\circ \psi_{23}\circ\psi_{12},
\qquad \psi\circ (u\tens 1) = 1\tens u,\quad \forall u\in P,
\end{equation}
\begin{equation}\label{diag.B}
(\id\tens\Delta)\circ\psi =
\psi_{12}\circ\psi_{23}\circ(\Delta\tens \id),
\qquad (\id\tens
\eps)\circ\psi =
\eps\tens \id.
\end{equation}
The triple $(P,C,\psi)$ is called an {\em entwining structure}.
\end{definition}
We will often use the notation $\psi(c\tens u) = u_\alpha\tens c^\alpha $
(summation over $\alpha$ is understood). In this notation
conditions (\ref{diag.A}) and (\ref{diag.B}) take a very simple
explicit form
$$
(uv)_\alpha \tens c^\alpha = u_\alpha v_\beta \tens
c^{\alpha\beta},
\qquad 1_\alpha\tens c^\alpha = 1\tens c,
$$
$$
u_\alpha \tens c^\alpha \sw 1\tens c^\alpha \sw 2 =
u_{\beta\alpha}\tens c\sw 1^\alpha\tens c\sw 2^\beta , \qquad
u_\alpha\eps(c^\alpha) = \eps(c) u,
$$
for any $u,v\in P$ and $c\in C$.

The entwining structure corresponds to an algebra factorisation in
the case $C$ finite-dimensional, built on $A=C^{*op}$ and $P$,
as explained in
\cite{BrzMa:coa}. Similarly, if $e\in C$ is grouplike, there is a
right coaction $\Delta_P:P\to P\tens C$ defined by
$\Delta_P(u)=\psi(e\tens u)$ and $M=M_e=\{u\in P|\
\Delta_P(u)=u\tens 1\}$ is a subalgebra. The map $\tilde\chi:P\tens
P\to P\tens C$ defined by $\tilde\chi(u\tens v)=u\Delta_P(v)$
descends to $\chi:P\tens_MP\to P\tens C$ and we have a {\em
copointed coalgebra bundle} $P(M,C,\psi,e)$ when $\chi$ is invertible and $P$.
This is the setting studied in
\cite{BrzMa:coa}.

We also note that for any entwining structure
we have a natural category $\M_P^C(\psi)$ of {\em
entwined modules}. The objects  are right $P$-modules and right
$C$-comodules $V$ such that for all $v\in V$, $u\in P$
 \begin{equation}\label{more.coa}
\Delta_P(v\ra u) = v\sw 0\ra\psi(v\sw 1\tens u) := v\sw 0\ra
u_\alpha \tens v\sw 1^\alpha,
\end{equation}
The morphisms are right $P$-module and right $C$-comodule maps.
The category $\M_P^C(\psi)$ generalises the category of unifying
or Doi-Koppinen modules \cite{Doi:uni}\cite{Kop:var} which
unifies various categories studied intensively in the
Hopf algebra theory (e.g. Drinfeld-Radford-Yetter (or crossed) modules, 
Hopf modules,
relative Hopf modules, Long modules etc.).

The algebra $P$
is an object in $\M_P^C(\psi)$, with the right regular action of
$P$ (by multiplication) if and only if there exists an
element $\tilde e\in P\tens C$ such that
\[ \tilde e\uo\psi(\tilde e\ut\tens \tilde e'\uo)\tens \tilde e'\ut
=(\id\tens\Delta)\tilde e,\quad (\id\tens\eps)\tilde e=1\]
(where $\tilde e'$ is another copy of $\tilde e$ and we use the notation $\tilde
e = \tilde e\su 1\tens \tilde e\su 2$, etc.). In this case
the coaction is
\[ \Delta_P(u)=\tilde e\uo\psi(\tilde e\ut\tens u),\quad\forall u\in P.\]
Notice that $\tilde e=\Delta_P(1)$. We then define\note{***}
\[ M=\{m\in P\,|\ \Delta_P(mv)=m\Delta_Pv\ \forall v\in P\}=\{m\in P\,
|\ \Delta_P(m)=m\tilde e\}\]
which is a subalgebra of $P$, and proceed as above, requiring $\chi$ to be
bijective. We will call this a general
{\em coalgebra bundle} $P(M,C,\psi)$. 
The copointed case corresponds to the
choice $\tilde e= 1\otimes e$.

There is also a converse: if $P$ is an algebra and a right
$C$-comodule, we say that the coaction is {\em Galois} if $M$
defined as above is such that $\chi$ is bijective. In this case
there is an entwining structure \cite{BrzHaj:coa}
\[ \psi(c\tens u)=\chi(c\uo\tens_M c\ut u),\quad \chi^{-1}(1\tens c)
=c\uo\tens c\ut,\quad \tilde e=\Delta_P(1)\]
and we have a coalgebra bundle. Because of these natural
properties, we will work now with these slightly more general
coalgebra bundles (or $C$-Galois extensions).  Our preliminary goal
in the present section is to make the evident generalisations of
the copointed theory in \cite{BrzMa:coa} to this case.

Next, a coalgebra bundle is {\em trivial} cf\cite{BrzMa:coa} (or
one says that the $C$-Galois extension is cleft) if there is a
convolution invertible map $\Phi:C\to P$ (the trivialisation or
cleaving map) such that
\begin{equation}
\Delta_P\circ\Phi = (\Phi\otimes \id)\circ\Delta.
\label{cov.phi}
\end{equation}
By considering the equality $1\sw 0\eps(c)\tens 1\sw 1 = 1\sw
0\psi(1\sw 1\tens\Phi(c\sw 1)\Phi^{-1}(c\sw 2))$ one finds that
\begin{equation}\label{cov.phi-1}
\psi(c\sw 1\tens \Phi^{-1}(c\sw 2)) = \Phi^{-1}(c)\Delta_P(1)
\end{equation}
which allows one to use the argument of the proof of
\cite[Proposition~2.9]{BrzMa:coa}  to show that $P\cong M\tens C$
as a left $M$-module and right $C$-comodule.

We turn now to the theory of connections, \note{*** moved defn OmegaP to
prelims} based on the theory for
the copointed case in \cite{BrzMa:coa}. As shown in
\cite[Proposition~2.2]{BrzMa:coa}, given an entwining structure
$(P,C,\psi)$ there is an entwining structure $(\Omega
P,C,\psi^\bullet)$, where
\[\psi^\bullet\mid_{C\otimes \Omega^{n-1}P} = \psi^n \equiv
\psi_{n,n+1}\psi_{n-1,n}\cdots\psi_{12}:C\tens P^{\tens n}\to
 P^{\tens n}\tens C\]
is the iterated entwining. Moreover,
\begin{equation}
\psi^\bullet\circ(\id\tens \extd) = (\extd\tens \id)\circ\psi^\bullet.
\label{cov.d}
\end{equation}
Therefore, given  $\tilde e:P\tens C$
we have $\Omega^{n}P\in \M_{\Omega P}^C(\psi)$
with the action right multiplication by $P$ and the coaction
\[\Delta_{\Omega^nP} = (\cdot_P\tens\id)(\id\tens\psi^{n+1})(\tilde e
\tens\id).\]

\begin{defin} A connection on $P(M,C,\psi)$ is a left $P$-module projection
$\Pi:\Omega^1P\to \Omega^1P$ such that (i) $\ker\Pi=P(\Omega^1M)P$
(ii) the map\note{***} $\Pi\circ \extd:P\to \Omega^1P$ commutes
with the right coaction.
\end{defin}

\begin{prop} Connections $\Pi$ are in 1-1 correspondence with
$\omega:C\to \Omega^1P$ such that

(i) $\tilde e\uo\omega(\tilde e\ut)=0$

(ii) $\tilde\chi\circ\omega(c)=1\tens c-\eps(c)\tilde e$

(iii) $\psi^2(c\o\tens \omega(c\t))=\omega(c\o)\tens c\t$.

The correspondence is via $\Pi(u\extd v)=u v\sw 0\omega(v\o)$ for
all $u,v\in P$.
\label{prop.connection.form}
\end{prop}
\proof Assume first that there is $\omega$ satisfying (i)-(iii). Then the
map $\Pi$ is well-defined since for all $u\in P$, $\Pi(u\extd 1) = u\omega
\tilde e\uo (\tilde e\ut) =0$, by (i). Next for any $u,v\in P$, $x\in M$
we have
\[
\Pi(u(\extd x) v)= \Pi(u\extd (xv)) - \Pi(ux\extd v) = u(xv)\sw 0 \omega
((xv)\sw 1) - uxv\sw 0\omega(v\sw 1) =0,
\]
since $\Delta_P$ is left $M$-linear. On the other hand, if $\sum_i u^i\extd v^i\in
\ker\Pi$, then using (ii) we have
\[
0 = \sum_i\tilde\chi(u^iv^i\sw 0\omega(v^i\sw 1)) = \sum_i(u^iv^i\sw 0 
\otimes v^i\sw 1
- u^iv^i\tilde e) = \sum_i\tilde\chi(u^i\extd v^i).
\]
Since $\ker\tilde\chi = P(\Omega^1M) P$, we have $\ker \Pi \subseteq P
(\Omega^1 M) P$,
i.e., $\ker\Pi = P(\Omega^1M)P$. Finally notice 
that for all $u\in P$,
$\Pi(\extd u) = u\sw 0\omega(u\sw 1)$. Therefore
\begin{eqnarray*}
\Delta_{\Omega^1P}(\Pi(\extd u)) & = & u\sw 0\psi^2(u\sw 1
\otimes \omega(u\sw 2)) \\
& = & u\sw 0\omega(u\sw 1)\otimes u\sw 2 \qquad
\qquad \mbox{\rm (by (iii))}\\
& = & \Pi(\extd u\sw 0)\otimes u\sw 1.
\end{eqnarray*}

Conversely, assume there is a connection in $P(M,C,\psi)$. This is
equivalent to the existence of a map $\sigma:P\otimes C^+\to \Omega^1 P$,
where $C^+=\ker\eps$, such that $\tilde\chi\circ\sigma = {\rm id}$ and 
$\Pi = \sigma\circ
\tilde\chi$. Define $\omega (c) = \sigma(1\otimes c - \eps(c)\tilde{e})$.
Clearly, (ii) holds. An immediate calculation verifies (i). The
definition of $\omega$ implies that $\Pi(u\extd v) = uv\sw 0\omega
(v\sw 1)$, for all $u,v\in P$.
Since $\Pi\circ\extd$ commutes with the coaction we have for all $u\in P$
\[
u\sw 0\psi^2(u\sw 1\otimes \omega(u\sw 2)) = u\sw 0\omega(u\sw 1)\otimes
u\sw 2.
\]
Since $\chi$ is bijective, for any $c\in C$ there is $c\su 1\otimes c\su 2
\in P\otimes_M P$ such that $c\su 1c\su 2\sw 0\otimes c\su 2\sw 1 = 1\otimes
c$. Thus we have
\begin{eqnarray*}
\psi^2(c\sw 1\otimes\omega(c\sw 2)) & = & c\su 1c\su 2\sw 0\psi^2(
c\su 2\sw 1\otimes \omega(c\su 2\sw 2))\\
& = & c\su 1c\su 2\sw 0\omega(c\su 2\sw 1)\otimes c\su 2\sw 2 =
\omega(c\sw 1)\otimes c\sw 2.
\end{eqnarray*}
Therefore $\omega$ satisfies (iii) and the proof of the proposition
is completed.
\endproof

Every connection $\Pi$, induces a {\em covariant
derivative}, $D = \extd -\Pi\circ\extd :P\to \Omega^1P$. In the
copointed case $D$ commutes with the right coaction, since $\extd$
itself commutes with the right coaction. 

\begin{prop} If $P(M,C,\psi,e)$ is a copointed trivial coalgebra bundle
with trivialisation $\Phi$ such that $\Phi(e)=1$,
\note{\tilde
e\uo\tens \Phi(\tilde e\ut)=1\tens 1$} 
and $\alpha:C\to \Omega^1M$ obeys $\alpha(e)=0$, 
\note{$\tilde
e\uo\alpha(\tilde e\ut)=0$}
 then
\[ \omega(c)=\Phi^{-1}(c\o)\alpha(c\t)\Phi(c\th)+\Phi^{-1}(c\o)\extd
\Phi(c\t)\]
is a connection.
\end{prop}
\proof
We verify directly that $\omega$ satisfies conditions (i)-(iii) of
Proposition~\ref{prop.connection.form} with $\tilde{e}=1\tens e$. 
\note{First notice that the normalisation condition $\Phi(e)=1$ implies that
$\Phi^{-1}(e) =1$.
\note{e\uo\tens \Phi(\tilde e\ut)=1\tens 1$ implies that $\tilde
e\uo\Phi^{-1}(\tilde e\ut\sw 1)\otimes \te\su 2\sw 2 =1\sw 0\otimes 1\sw
1 = \te$. }}
We have
\[\omega(e) =  \Phi^{-1}(e)\alpha(e)\Phi(e)+\Phi^{-1}(e)\extd
\Phi(e) = \extd 1  = 0.\]
Next, take any $c\in C$ and compute
\begin{eqnarray*}
\tilde\chi\circ\omega(c) & = & \tilde\chi(\Phi^{-1}(c\sw 1)\tens
\Phi(c\sw 2)-\eps(c)1\otimes 1)\\
& = &
\Phi^{-1}(c\sw 1)\Phi(c\sw 2)\tens c\sw 3 - \eps(c) 1\tens e
= 1\tens c - \eps(c)\te,
\end{eqnarray*}
where we used that the first summand in $\omega$ is in $P(\Omega^1M)P$.
Finally we have
\begin{eqnarray*}
\psi^2(c\sw 1\otimes\omega(c\sw 2)) &=& \psi^2(c\sw 1\otimes
\Phi^{-1}(c\sw 2)\alpha(c\sw 3)\Phi(c\sw 4))+ \psi^2(c\sw 1\otimes
\Phi^{-1}(c\sw 2)\extd \Phi(c\sw 3)) \\
& = & \Phi^{-1}(c\sw 1)\psi^2(e\tens\alpha(c\sw 2)\Phi(c\sw
3)) +
\Phi^{-1}(c\sw 1)\psi^2(e\tens \extd \Phi(c\sw 2)) \\
& = & \Phi^{-1}(c\sw 1)\alpha(c\sw 2)\Delta_P(\Phi(c\sw
3))+\Phi^{-1}(c\sw 1)\extd\Phi(c\sw 2)\tens c\sw 3 \\
& = & \omega(c\sw 1)\tens c\sw 2,
\end{eqnarray*}
where we used that $\Omega^1P\in \M_{\Omega P}^C(\psi\sp\bullet)$ and (\ref{cov.phi-1}) to
derive the second equality, and that $\alpha(c)\in \Omega^1M$, $\Phi$ is
an intertwiner and (\ref{cov.d}) to derive the third one.
 \endproof

For another class of examples one has coalgebra homogeneous spaces
associated to coalgebra surjections $\pi:P\to C$.
Thus, let $P$ be a Hopf algebra and $M$ a subalgebra of $P$ such
that $\Delta(M)\subset P\tens M$ (an embeddable $P$-homogeneous
quantum space). Define the quotient coalgebra $C = P/(M^+P)$, where
$M^+ = \ker\eps\cap M$ is the augmentation ideal. There is a
natural right coaction of $C$ on $P$ given as $\Delta_P= (\id \tens
\pi)\circ
\Delta$, where $\pi: P\to C$ is the canonical surjection.
It is clear that $M\subseteq \{u\in P| \Delta_Pu=u\tens e\}$, with
$e=\pi(1)$, and we assume that this is an equality (this is known
to hold for example if\cite{Sch:nor}  $P$ is faithfully flat as a
left $M$-module). Then $P(M,C,\pi(1))$ is a coalgebra bundle.
Since $\tilde e=1\tens\pi(1)$ we have $e=\pi(1)$,
i.e. a copointed coalgebra bundle as in \cite{BrzMa:coa}. In this
case we know that if $i:C\to P$ is a linear splitting of $\pi$ such
that
\[ i(e)=1,\quad \eps\circ i=\eps,
\quad i(c\t)\t\tens c\o\ra(Si(c\t)\o)i(c\t)\th)=i(c\o)\tens c\t\]
then
\[ \omega(c)=Si(c)\o\extd i(c)\t\]
is a left-invariant connection and every left-invariant connection
on the bundle is of this form (cf.\ \cite{wor}). 
The left-invariance here means that
$\Delta_{\Omega^1P}\omega(c) = 1\otimes \omega(c)$ for all $c\in C$. 
We use here the right action of $P$
on $C$ given by $c\ra u=\pi(vu)$ for any $v\in\pi^{-1}(c)$.

The theory of connections can be developed also for nonuniversal calculi
$\Omega^1(P)=\Omega^1P/\CN$ where $\CN\subset\Omega^1P$ is a
sub-bimodule, although the situation is slightly more complicated. 
We say that $\Omega^1(P)$ is a {\em differential calculus on
$P(M,C,\psi)$} iff it is {\em covariant} in
the sense
\[ \psi^2(C\tens \CN)\subseteq \CN\tens C\]
so that the coaction $\Delta_{P\tens P}$ descends to $\Omega^1(P)$.
This is obtained from $\psi^2_{\CN}$ defined by
\[ \psi^2_{\CN}\circ (\id\tens \pi_{\CN})=(\pi_{\CN}\tens\id)\circ
\psi^2\]
where $\pi_{\CN}:\Omega^1P\to \Omega^1(P)$ is the canonical
surjection. We have
\[ \Delta_{\Omega^1(P)}=\tilde e\uo\psi^2_{\CN}(\tilde e\ut\tens(\ )).\]
Let $\CM=(P\tens C^+)/\tilde\chi(\CN)$ (and denote by $\pi_{\CM}$
the canonical surjection). This is a left $P$-module (since
$\tilde\chi$ is left $P$-module map) by $u\la
m=\sum_i\pi_{\CM}(uv_i\tens c_i)$ for any $\sum v_i\tens
c_i\in\pi^{-1}_{\CM}(m)$. We can then define
\[ \Lambda=\{\lambda\in\CM|\exists c\in C,\ {\rm s.t.}\ \lambda
=\pi_{\CM}(1\tens c
-\eps(c)\tilde e)\}.\]
The action provides a surjection $P\tens\Lambda\to\CM$.

\begin{defin} A connection with a nonuniversal calculus is a
left $P$-module projection
$\Pi:\Omega^1(P)\to \Omega^1(P)$ such that
$\ker\Pi=\Omega^1(P)_{\rm hor}$ and $\Pi\circ\extd$ 
commutes with the right coaction.
\end{defin}

Here $\Omega^1(P)_{\rm hor}=P(\extd M)P$. As usual, we define
$\chi_{\CN}\circ\pi_{\CN}=\pi_{\CM}\circ\tilde \chi$. It is a left $P$-module
map and the sequence
\[ 0\to \Omega^1(P)_{\rm hor}\to \Omega^1(P){\buildrel 
\chi_{\CN}\over\to}\CM\to 0\]
is exact.

\begin{prop} Suppose $P\tens\Lambda\cong\CM$ by the surjection above. Then
 connections $\Pi$ on $\Omega^1(P)$ are in 1-1 correspondence with
$\omega:\Lambda\to \Omega^1(P)$ such that

(i) $\chi_{\CN}\circ\omega=1\tens \id$

(ii) $\psi^2_{\CN}(c\o\tens
\omega(\pi_\Lambda (c\t)))=\omega(\pi_\Lambda(c\o))\tens c\t$ where
$\pi_\Lambda(c)= \pi_{\CM}(1\tens c-\eps(c)\tilde e)$.

The correspondence is via $\Pi(u\extd v)=u \sum_i v_i
\omega(\lambda_i)$ for all $u,v\in P$ and $\sum_i
v_i\tens\lambda_i\in P\tens\Lambda$ such that $\sum_i
v_i\la\lambda_i=\chi_{\CN}(\extd v)$.
\end{prop}
\proof The proof is analogous to the proof of
Proposition~\ref{prop.connection.form}. \endproof

In the case of a homogeneous bundle where $P$ is a Hopf algebra and
$e=\pi(1)$, a natural type of calculus $\Omega^1(P)$ is a
left-covariant one defined by an ideal $\CQ$ in
$\ker\eps\subset P$.

\begin{example} For a homogeneous bundle with left-covariant calculus, 
$\Lambda=C^+/\pi(\CQ)$ and $P\tens\Lambda\cong\CM$. Moreover if for all
$q\in \CQ$, $u\in P$, $q\sw 2\tens \pi(u(Sq\sw 1)q\sw 3)\in \CQ\tens C$,
then $\Omega^1(P)$ is a calculus on $P(M,C,\psi,\pi(1))$. In particular,
if $\Omega^1(P)$ is a bicovariant calculus on $P$ then it is a 
calculus on $P(M,C,\psi,\pi(1))$. 
\note{ Moreover, let  $i:\Lambda\to P^+/\CQ$ be a
splitting of the projection $\bar\pi:P^+/\CQ\to\Lambda$ induced by $\pi$
and such that
\[ \pi_{\CQ}(\iota(\pi_\Lambda(c\t))\t)\tens c\o\ra \iota
(\pi_\Lambda(c\t))\o=  i(\pi_\Lambda(c\o))\tens c\t,\]
where $\pi_{\CQ} :  P\to P/\CQ$ is the canonical surjection and
$\iota(\lambda)\in \pi^{-1}(i(\lambda))$. Then
\[ \omega(\lambda)=S\iota(\lambda)\o\extd \iota(\lambda)\t\]
is a left-invariant connection. Every left-invariant connection
is of this form.}
\end{example}
\proof Recall that any element $n\in \CN$ is of the form 
$n = \sum_iu^iSq^i\sw
1\otimes q^i\sw 2$ for some $u^i\in P$, $q^i\in \CQ$. For any $u\in P$,
$q\in \CQ$ we have $u\tens \pi(q) = \tilde\chi(uSq\sw 1\tens q\sw 2)\in
\tilde\chi(\CN)$. On the other hand $\tilde\chi(\sum_iu^iSq^i\sw
1\otimes q^i\sw 2) = \sum_i u^i\tens \pi(q^i)\in P\tens \CQ$. This
proves that $\tilde\chi(\CN) = P\tens \pi(\CQ)$. Therefore $\CM = P\tens
C^+/\tilde\chi(\CN) = P\tens (C^+/\pi(\CQ))$, and $\Lambda =
C^+/\pi(\CQ)$.
\note{we
can consider the following commutative diagram:
$$
\begin{CD}
@. 0 @>>>    0 @>>>  0 @>>> 0\\
  @.      @AAA   @AAA @AAA\\
0 @>>>    \tilde{\chi}(\CN )   @>>>    P\tens C^+ @>>> \CM @>>> 0 \\
  @.      @|   @| @AA{s}A\\
0 @>>> P\tens\pi(\CQ)    @>>> P\tens C^+ @>>> P\tens
  \Lambda @>>> 0\\
@.     @AAA   @AAA @AAA\\
0@>>> 0 @>>>    0    @>>> \ker s @.
\end{CD}
$$
The first row is the definition of $\CM$, while the second one is the
definition of $\Lambda$. 
The application of the Snake Lemma (cf. \cite[Section~1.2]{Bou:com})
to the first of the above diagram
yields $\ker s =0$, so that $\CM$ is isomorphic to $P\otimes\Lambda$.}

Finally, take any $c\in C$
and let $v\in \pi^{-1}(c)$. We have:
\begin{eqnarray*}
\sum_i\psi^2(c\otimes u^iSq^i\sw 1 \!\!\!\!&\otimes &\!\!\!\! q^i\sw
2) =  \sum_i u^i\sw
1Sq^i\sw 2\tens \psi(\pi(vu^i\sw 2 Sq^i\sw 1)\tens q^i\sw 3)\\
& = & \sum_i u^i\sw 1Sq^i\sw 2\tens q^i\sw 3\tens \pi(vu^i\sw 2(Sq^i\sw
1)q^i\sw 4).
\end{eqnarray*}
By the assumption on $\CQ$ the last expression is in $\CN\tens \CQ$, so
that the resulting calculus $\Omega^1(P)$ is a calculus on
$P(M,C,\psi,\pi(1))$. If $\CQ$ defines a bicovariant calculus then $\CQ$
is Ad-stable, so that the required condition is immediately satisfied.
\endproof

\section*{\normalsize\sc\centering 4. Bijectivity of $\psi$ and strong
connections}\vspace{-.6\baselineskip}
\setcounter{section}{4}
\setcounter{definition}{0}

In this section we return to some technical considerations. 
For simplicity here and
in most of what follows, we will concentrate on the universal
differential calculus. First of all, we consider the question
of when $\psi$ is bijective. It
plays the role in the Hopf algebra case of having a bijective antipode,
and
allows us to relate left and right handed versions of the theory.

\begin{lem} If $\psi$ is bijective then $P$ is a left $C$-comodule by
\[ {}_P\Delta(u)=\psi^{-1}(u\tilde e).\]
Moreover, $M=\{u\in P|\ {}_P\Delta u=\psi^{-1}(\tilde e) u\}$.
\end{lem}
\proof This lemma is part of \cite[Lemma~6.5]{Brz:mod}. \endproof

In the copointed case, it is easy to see that if $\psi$ is
bijective then $P^{\tens(n+1)}$ is a left $C$-comodule by ${}_{P^{\tens
(n+1)}}\Delta=\psi^{-(n+1)}((\ )\tens e)$. This coaction restricts
to $P\tens M^{\tens n}$ and $\Omega^nP$.

\begin{prop} In the copointed case, let
$\omega$ be a connection on $\Omega^1P$ with $\psi$ bijective. Then
$\bar\Pi:\Omega^1P\to \Omega^1P$ defined by
\[ \bar\Pi((\extd u)v)=\omega(u\sw{1})u\sw{\infty} v\]
is a right-connection in the sense

(i) $\bar D=(\id-\bar\Pi)\circ\extd$ is a left $C$-comodule map.

(ii) $\bar\Pi$ is a right $P$-module projection and
$\ker\bar\Pi=P(\Omega^1M)P$.

\end{prop}
\proof (i) We introduce the notation $\psi^{-1}(u\tens c) =
c_\alpha\tens u^\alpha$, for all $c\in C$, $u\in P$. One easily finds
that
\begin{equation}
c_\alpha\sw 1\tens c_\alpha\sw 2 \tens u^\alpha = c\sw 1_\alpha\tens
c\sw 2_\beta\tens u^{\alpha\beta}
\label{psi-1.a}
\end{equation}
 and
${}_P\Delta(u) = e_\alpha\tens u^\alpha$. We have
\begin{eqnarray*}
\psi^2(u\sw 1\otimes\omega(u\sw 2)u\sw\infty) & = & \omega(u\sw
1)\psi(u\sw 2\otimes u\sw\infty) \\
& = & \omega(e_\alpha\sw 1)\psi(e_\alpha\sw 2\tens u^\alpha) \\
& = & \omega(e_\alpha)\psi(e_\beta\tens u^{\alpha\beta})\\
& = & \omega(e_\alpha)u^\alpha\tens e = \omega(u\sw 1)u\sw\infty\tens e,
\end{eqnarray*}
where we used that $e$ is group-like and (\ref{psi-1.a}) to derive the
third equality. This implies that
\[\psi^{-2}(\omega(u\sw 1)u\sw \infty \tens e) = u\sw 1\tens \omega(u\sw
2)u\sw\infty ,\]
which is precisely the left $C$-covariance of $\bar\Pi\circ\extd$ and,
consequently, implies the left-covariance of $\bar D$.

(ii) It is clear that $\bar\Pi$ is a right $P$-module map. The following
diagram commutes:
\[
\begin{CD}
0@>>> P(\Omega^1M)P @>>> \Omega^1P @>{\tilde\chi}>> P\otimes C^+@>>> 0\\
@. @VV{=}V @VV{=}V @VV{\psi}V \\
0@>>> P(\Omega^1M)P @>>> \Omega^1P @>{\tilde\chi_L}>> C^+\otimes P @>>> 0
\end{CD}
\]
where $\tilde\chi_L = \psi^{-1}\circ\tilde\chi$ (explicitly,
$\tilde\chi_L(u\tens v) = u\sw 1\tens u\sw\infty v$). Since $P$ is
a coalgebra principal bundle the top sequence is exact. Furthermore
 $\psi$ is bijective and ${}_P\Delta$ is right $M$-linear thus the
bottom sequence is also exact. It is
split by the map $\sigma: C^+\otimes P\to \Omega^1P$,
$\sigma(c\otimes u) = \omega(c) u$. Indeed,
\[\tilde\chi_L\circ\sigma(c\otimes u) = \tilde\chi_L(\omega(c))u =
\psi^{-1}(1\otimes c)u = c\otimes u,
\]
where we used that $\tilde\chi_L$ is a right $P$-module map and that $\omega$
is a connection one-form (Proposition~\ref{prop.connection.form}(ii)).
Now notice that $\bar\Pi = \sigma\circ \tilde\chi_L$, and the fact that
$\sigma$ is a splitting (i.e. $\tilde\chi_L\circ\sigma = \id$) of the
above sequence implies both that $\bar\Pi$ is a projection and has the
kernel as stated.
\endproof

Finally, a connection is {\em strong} if $(\id-\Pi)\circ\extd$ has its
image in $(\Omega^1M)P$ \cite[Definition~2.1]{Haj:str}. 
These are the connections most closely associated to
the base and used in the theory of associated bundles etc. Recently,
a simple condition for strongness was given in the Hopf
algebra case, in \cite{Ma:rie}. This can be generalised to the coalgebra
case. 

\begin{prop} A connection  on a copointed coalgebra
bundle $P(M,C,\psi, e)$ is  strong
{\em iff}
\begin{equation}
(\id\tens\Delta_P)\omega(c)=
1\tens 1\tens c-\eps(c)1\tens 1\tens e+\omega(c\o)\tens c\t.
\label{l.strong.1}
\end{equation}
Furthermore, if $\psi$ is bijective then a connection
is strong {\em iff}
\[ ({}_P\Delta\tens\id)\omega(c)=c\tens 1\tens 1-e\tens 1\tens 1\eps(c)
+c\o\tens\omega(c\t).\]
\label{l.strong.l}
\label{l.strong.r}
\end{prop}
\proof Assume that $\omega$ is strong. This is equivalent to the
statement  that
\begin{equation}
(\id\tens \Delta_P)\circ D(u) = \Delta_{\Omega^1P}\circ D(u), \qquad   
\forall u\in P.
\label{l.strong.2}
\end{equation}
Using the explicit
definition of $\extd$ and $D$, Proposition~\ref{prop.connection.form}(iii), 
 as well as the fact that $\Omega^1P\in
\M_{\Omega P}^C(\psi^\bullet)$  one finds that
(\ref{l.strong.2}) implies that
\[
(\id\tens\Delta_P)(u\sw 0\omega(u\sw 1)) = u\sw 0\tens 1\tens u\sw
1-u\tens 1\tens e +u\sw 0\omega(u\sw 1)\tens u\sw 2.\]
Next for all $c$, let $c^\su 1\tens c\su 2\in P\otimes_MP$ be the
translation map, i.e. $c^\su 1\tens c\su 2 = \chi^{-1}(1\otimes c)$. It
means that
$c\su 1c\su 2\sw 0\otimes c\su2\sw 1 = 1\otimes c$. Using the above
equality and the fact that $c\su 1c\su 2 = \eps(c)$, we have
\begin{eqnarray*}
(\id\tens\Delta_P)\circ\omega(c) & = & (\id\tens\Delta_P)(c\su 1c\su 2\sw
0\omega(c\su 2\sw 1)) = c\su 1(\id\tens\Delta_P)(c\su 2\sw
0\omega(c\su 2\sw 1))\\
& = & c\su 1c\su 2\sw 0\tens 1\tens c\su 2\sw
1-c\su 1c\su 2\tens 1\tens e \\
&&+c\su 1c\su 2\sw 0\omega(c\su 2\sw 1)\tens
 c\su 2\sw 2 \\
& = & 1\tens 1\tens c-\eps(c)1\tens 1\tens e+\omega(c\o)\tens c\t ,
\end{eqnarray*}
i.e. (\ref{l.strong.1}) holds. Conversely, an easy calculation reveals
that (\ref{l.strong.1})  implies (\ref{l.strong.2}), i.e.,
the connection is strong as required.

 The second assertion is obtained by applying $\psi^{-2}$ to
(\ref{l.strong.1}). \endproof

As in \cite{Ma:rie}, the significance of this is that this is manifestly
a `strongness' condition for the left-handed theory with $\bar\Pi$.
In studying the coalgebra frame resolutions we will need both the left and the
right handed theories simultaneously, and we see that if one holds so does
the other for a given $\omega$.

A situation where $\psi$ is bijective is a homogeneous bundle $\pi:P\to C$
with $P$ having bijective antipode.

 \begin{prop} For a homogeneous coalgebra bundle with bijective
antipode, strong left-invariant connections are in 1-1
correspondence with splittings $i:C\to P$ of
$\pi$ which are covariant with respect to $(\id\tens\pi)\circ\Delta$
and $(\pi\tens\id)\circ\Delta$, and such that $i(\pi(1))=1$ and 
$\eps\circ i=\eps$. In this
case
\[ \omega(c)=Si(c)\o\extd i(c)\t.\]
\end{prop}
\proof Given such a splitting $i:C\to P$ of $\pi$, consider
$\omega(c)= Si(c)\o\extd i(c)\t$
as stated. The normalisation conditions imply that
 $\omega(\pi(1)) = 0$ and
$\tilde{\chi}\circ\omega(c) = 1\otimes c -
\eps(c)1\otimes\pi(1)$. Also
\begin{eqnarray*}
\psi^{2}(c\sw 1\otimes \omega(c\sw 2)) & = & Si(c\sw 2)\sw 2\extd
i(c\sw 2)\sw 3\tens\pi(i(c\sw 1) Si(c\sw 2)\sw 1i(c\sw 2)\sw 4)\\
& = &
Si(c)\sw 3\extd
i(c)\sw 4\tens\pi(i(c)\sw 1Si(c)\sw 2i(c)\sw 5) \quad \mbox{\rm ($i$
is left-covariant)}\\
& = & Si(c)\o\extd i(c)\t\tens \pi(i(c)\sw 3)\\
& = & Si(c\sw 1)\sw 1\extd i(c\sw 1)\sw 2\tens \pi(i(c\sw 2)) \qquad
\mbox{\rm ($i$ is right-covariant)}\\
& = & \omega(c\sw 1)\tens c\sw 2 \qquad \qquad\qquad\qquad\qquad
\mbox{\rm ($\pi$ is split by $i$)}
\end{eqnarray*}
Proposition~\ref{prop.connection.form} implies that $\omega$ is a
connection  one-form. Finally, compute
\begin{eqnarray*}
(\id\otimes\Delta_P)(\omega(c)) & = & Si(c)\sw 1\tens i(c)\sw
2\tens\pi(i(c)\sw 3) - \eps(c) 1\tens 1\tens \pi(1)\\
& = & Si(c\sw 1)\sw 1\tens i(c\sw 1)\sw 2\tens c\sw 2 - \eps(c) 1\tens
1\tens\pi(1)\\
& = & \omega(c\sw 1)\tens c\sw 2+1\tens 1\tens c - \eps(c) 1\tens 1\tens
\pi(1),
\end{eqnarray*}
where the use of the fact that $i$ is a right covariant splitting was
made in the derivation of the second equality.
Proposition~\ref{l.strong.r} now implies that the connection
corresponding to $\omega$ is strong.

Conversely, assume that there is a  strong connection with the
left-invariant connection form $\omega$. Then  the left-invariance
of $\omega$ implies that there exists a splitting $i: C\to P$ of $\pi$
such that $\eps\circ i =\eps$ and
$\omega(c)=Si(c)\o\extd i(c)\t$  (cf.
\cite[Proposition~3.5]{BrzMa:dif}). The fact that $\omega(\pi(1)) =0$
implies that $i(\pi(1))=1$. Applying $(\id\tens \Delta_P)$ to
this $\omega$ and using Proposition~\ref{l.strong.r} one deduces that
$i$ is right-covariant. Bijectivity of $S$ implies that $\psi$ is
bijective (cf. \cite{Brz:mod}). The left coaction induced by
$\psi^{-1}$ is ${}_P\Delta(u) = \pi(S^{-1}u\sw 2)\tens u\sw 1$. By
Proposition~\ref{l.strong.r}
\[ ({}_P\Delta\tens P)\omega(c) = \pi(i(c)\sw 1) \tens Si(c)\sw 2\tens i(c)\sw
3 -\eps(c) \pi(1)\tens 1\tens 1
\]
must be equal to
\[
c\sw 1\tens Si(c\sw 2)\sw 1\tens i(c\sw 2)\sw 2 -\eps(c)\pi(1)\tens
1\tens 1.
\]
Applying $\id\tens S^{-1}\tens \eps$ to this equality one deduces that
$i$ must be left-covariant. This completes the proof.
\endproof

This is the analogue for coalgebra bundles of the bicovariant formulation of
strong canonical connections in the Hopf algebra case in \cite{HajMa:pro}.

\section*{\normalsize\sc\centering 5. Frame resolutions, covariant
derivatives and torsion}\vspace{-.6\baselineskip}
\setcounter{section}{5}
\setcounter{definition}{0}

In this section we define frame resolutions 
in the coalgebra setting, following the theory
introduced recently in \cite{Ma:rie} in the Hopf algebra case. The
theory depends heavily on the notion of associated bundles, so we
recall these briefly. In the coalgebra case there are two kinds of
associated bundles (which are equivalent in the Hopf algebra case),
as studied recently in \cite{Brz:mod}.

\bde\label{vc-ext}
Let $P(M,C)$ be a coalgebra bundle.

(i) The left associated bundle (or module) to a left $C$-comodule $V$
is $E=P\Box_C V$.

(ii) The right associated bundle (or module) to a right $C$-comodule
$V$ is $\bar E=(V\tens P)_0$, the fixed
subobject, where $V\tens P$ is an object of $\M_P^C(\psi)$ by
multiplication from the right and $\Delta_{V\tens P}(v\tens u)=v\sw
0\tens\psi(v\o\tens u)$.
\ede

The cotensor product $W\Box_{C}V$ here, between a left comodule $V$ and
right comodule $W$ is defined by the exact sequence\cite{MilMoo:str}
\[
0\longrightarrow W\Box_{C}V\hookrightarrow W\ot V\st{\Delta_W\ot 
\id-\id\ot\, {}_V\Delta}
{\longrightarrow}W\ot C\ot V.\]
This is just the arrow reversal of the usual
tensor product. Less conventional is the fixed subobject
\[(V\tens P)_0 =\{\sum_i v_i\tens u_i\in V\tens P|\ v_i\sw
0\tens \psi(v_i\o\tens u_i)=v_i\tens u_i\tilde e\uo\tens\tilde
e\ut\}.\] This is the natural analogue for
coalgebra bundles of the associated bundles in the quantum group
gauge theory of \cite{BrzMa:gau}.

\ble\label{lem.inv}
For a copointed coalgebra bundle $P(M,C,\psi,e)$, let $(P\otimes
M^{\otimes n})_0=\{w\in P\otimes M^{\otimes n}|\ \psi^{n+1}(e\tens
w) = w\tens e\}$ be the {\em invariant} subset of $P\otimes
M^{\otimes n}$. If $\psi$ is bijective then $(P\otimes M^{\otimes
n})_0
= M^{\otimes{n+1}}$.
\ele
\proof Clearly $M^{\otimes{n+1}}\subset (P\otimes M^{\otimes n})_0$. If
$w\in (P\otimes M^{\otimes n})_0$ then $\psi^{n+1}(e\tens w) =
w\tens e$. Applying $\psi^{-(n+1)}$ one deduces that
$\psi^{-(n+1)}(w\tens e)
= e\tens w$. Let $w = \sum_i u^i\otimes m_1^i\otimes \cdots \otimes
m_n^i$. Since for all $m\in M$, $\psi^{-1}(m\otimes e) = e\otimes m$
one immediately finds that $e\otimes \sum_i u^i\otimes m_1^i\otimes
\cdots \otimes m_n^i = \sum_i \psi^{-1}(u^i\otimes e) \otimes
m_1^i\otimes \cdots \otimes m_n^i$. This in turn implies that for
all $i$, $u^i\in M$. \endproof

Now we can extend the notion of a strongly horizontal
form from \cite{BrzMa:gau}
\bde
Let $E$ be a left bundle associated to a  copointed coalgebra bundle
$P(M,C,\psi,e)$ and a left $C$-comodule $V$.
A {\em right strongly tensorial n-form} on $E$ is a
linear map $\phi : V\to P(\Omega^nM)$ such that
\begin{equation}\label{con.phi.n}
\psi^{n+1}\circ(\id\tens\phi)\circ{}_V\Delta = \phi\tens e,
\end{equation}
\ede
By the extension of the notation above, the space of right strongly 
tensorial $n$-forms will be denoted by
${\rm Hom}_0(V,P(\Omega^n M))$ (in \cite{Brz:mod} right strongly
0-forms ${\rm Hom}_0(V,P)$ are denoted by ${\rm Hom}_\psi(V,P)$).
${\rm
Hom}_0(V,P(\Omega^n M))$ has a right $M$-module structure
defined by $(\phi\cdot m)(v) = \phi(v)m$.

\begin{prop}\label{pro.str} Let $P(M,C,\psi,e)$ be a copointed
coalgebra bundle with
$\psi$ bijective and $P$ flat as a right $M$-module (or $V$-coflat as a
left $C$-comodule).
Then right strongly tensorial forms ${\rm Hom}_0(V,P(\Omega^n M))$
and  ${}_M\Hom(E,\Omega^nM)$ are isomorphic 
 as right $M$-modules.
\end{prop}
\noindent
\proof The proof of this proposition is analogous to the proof of
\cite[Theorem~4.3]{Brz:mod}. We include it here for completeness.
The flatness (coflatness) assumption implies that
$(P\otimes_MP)\square_CV \cong
P\otimes_M(P\square_CV)$, canonically (cf. \cite[p.~172]{Sch:pri}).
Thus there is a left
$P$-module isomorphism $\rho: P\otimes_M E \to P\otimes V$,
obtained as a composition of $\chi\otimes\id$ with the canonical
isomorphism $P\otimes C\square_CV \stackrel{\sim}{\to} P\otimes V$, i.e.,
$\rho = \cdot\otimes \id$, $\rho^{-1} = (\chi^{-1}\otimes
\id)\circ(\id\otimes{}_V\Delta)$.
Following \cite{DoiTak:miy}, apply  ${\rm Hom}_{P}(-,P(\Omega^n M))$ to
$\rho$ to deduce the right $M$-module isomorphism
$
{\rm Hom}(V,P(\Omega^n M)) \stackrel{\sim}{\to}
{\rm Hom}_{M}(E,P(\Omega^nM))$, given by $\phi\mapsto s_{\phi}$, 
$s_{\phi} (\sum_iu^i\tens v^i)=
\sum_iu^i\phi(v^i)$.
For any $\phi\in {\rm Hom}(V,P(\Omega^nM))$, $x=\sum_iu^i\tens v^i \in
E$ we have
\begin{eqnarray*}
\Delta_{\Omega P}(s_{\phi}(x)) & = &\sum_i\Delta_{\Omega P}
(u^i\phi(v^i))=
\sum_i u^i\sw 0\psi^{n+1}(u^i\sw 1\tens\phi(v^i)) \\
& = &\sum_i u^i\psi^{n+1}(v^i\sw{1}\tens\phi(v^i\sw{\infty})),
\end{eqnarray*}
since $\sum_iu^i\sw 0\tens u^1\sw 1\tens v^i = \sum_iu^i\tens v^i\sw
1\tens v^i\sw\infty$ by the definition of $E=P\Box_CV$. 
By Lemma~\ref{lem.inv}, $\Omega^nM = (P(\Omega^nM))_0$, therefore
 $s_{\phi}(x)\in \Omega^nM$ iff
\begin{equation}
\sum_i u^i\psi^{n+1}(v^i\sw{1}\tens\phi(v^i\sw{\infty})) = \sum_i
u^i\phi(v^i)\tens e.
\label{sec.1}
\end{equation}
Clearly, (\ref{con.phi.n}) implies (\ref{sec.1}). Applying (\ref{sec.1})
to $\rho^{-1}(1\otimes v)$
one easily finds that
(\ref{sec.1}) implies (\ref{con.phi.n}). Therefore the right $M$-module 
isomorphism
$\phi\mapsto s_{\phi}$ restricts to the isomorphism
${\rm Hom}_0(V,P(\Omega^nM))\stackrel{\sim}{\to}
{\rm Hom}_{M}(E,\Omega^nM)$ as required.
\endproof

\note{\bre
{\rm It can be easily shown that ${\rm Hom}_0(V,A(\Omega^n B))$ is a
right $B$-submodule of the space ${\rm Hom}_{\psi^{n+1}}(V,A(\Omega^n
B))$ of linear maps $\phi:V\to A(\Omega^n B)$ characterised by the
condition
\begin{equation}
\psi^{n+1}\circ(C\tens\phi)\circ{}_V\Delta = \phi\tens \eta_e. \label{lstr}
\end{equation}
\note{ and  $\psi^{n+1} = (A^{\otimes n}\tens\psi)\circ(A^{\otimes
n-1}\tens\psi\tens
A)\circ\ldots\circ(\psi\tens A^{\otimes n})$.}
\sloppypar Using the same
techniques as in the proof of Proposition~\ref{pro.str} one easily
finds that module ${\rm Hom}_{\psi^{n+1}}(V,A(\Omega^n
B))$ is isomorphic to a right $B$-module ${\rm
Hom}_{B-}(E,\overline{A(\Omega^nB)})$, where $\overline{A(\Omega^nB)}$
is the {\em invariant} subset of $A(\Omega^nB)$, {\em i.e.} the set of
all $w\in A(\Omega^nB)$ such that $\psi^{n+1}(e\tens w) = w\tens
e$. Clearly $\Omega^nB\subseteq \overline{A(\Omega^nB)}$. If this
inclusion is an equality then all strongly tensorial forms are fully
characterised by (\ref{lstr}).}
\ere}

Proposition~\ref{pro.str} is the coalgebra bundle version of
\cite[Lemma~3.1]{Ma:rie}, and allows us to define similarly,

\bde\label{frame.resolution} cf\cite[Definition~3.2]{Ma:rie}
A {\em coalgebra frame resolution} of an algebra $M$ is a left
bundle  $E$ associated to a copointed coalgebra bundle $P(M,C,\psi,e)$
with bijective $\psi$, and $V$,
together with a right strongly tensorial one-form $\theta: V\to
P(\Omega^1M)$ such that $s_\theta:E\to\Omega^1M$ corresponding
under Proposition~4.4 is an
isomorphism of left $M$-modules.
\ede

As in \cite{Ma:rie}, we can now
proceed to deduce the left $M$-module isomorphism
\begin{equation}
\id\otimes s_\theta : (\Omega^1M)P\Box_CV \cong \Omega^1M\otimes_M
\Omega^1M =\Omega^2M .\label{isom}
\end{equation}
Here, the cotensor product is defined with respect to the right
coaction $\Delta_{(\Omega^1M)P}:(\Omega^1M)P\to (\Omega^1M)P\tens C$
given by $\Delta_{(\Omega^1M)P}(w) = \psi^2(e\tens w)$ (it is an
easy exercise which uses (\ref{diag.A}) to verify that
$(\Omega^1M)P$ is closed under this coaction).

Furthermore, given a frame resolution, we can now define a {\em
covariant derivative} $\nabla
:\Omega^1M\to
\Omega^2M$ corresponding to a strong connection $\Pi$ in 
$P(M,C,\psi,e)$ by \cite[Proposition~3.3]{Ma:rie}
\begin{equation}\label{cov.der}
\nabla = (\id\tens s_\theta)\circ (D\Box_C\id)\circ s_\theta^{-1} 
:\Omega^1M\to
\Omega^2M.
\end{equation}
The map $\nabla$ is well-defined since
$D$ is an intertwiner so that the expression
$D\Box_C\id$ makes sense. Furthermore, by the strongness
assumption $D(P)\subseteq (\Omega^1M)P$ so the isomorphism
(\ref{isom}) implies that the output of $\nabla$ is in $\Omega^2M$.
Finally, it can be easily verified (cf.
\cite[Proposition~3.3]{Ma:rie}) that $\nabla(m\cdot w) =
m\cdot\nabla w + \extd m\otimes_M w$, for any $m\in M$ and $w\in
\Omega^1M$, so that $\nabla$ is a connection on $\Omega^1M$ as a left 
$M$-module.

Next, cf \cite[Proposition~3.5]{Ma:rie}, we define the torsion of
a connection $\nabla$ by
\[T
=\extd - \nabla: \Omega^1M\to \Omega^2M.\]
By Proposition~\ref{pro.str} this $T$ can be also viewed as a map
$\T :V\to P(\Omega^2M)$ provided $P$ is $M$-flat.

\begin{prop} If $\omega$ is a strong connection on $P(M,C,\psi,e)$ and
$\psi$ is bijective then there is a covariant derivative
\[ \bar D:\Hom_0(V,P\Omega^nM)\to \Hom_0(V,P\Omega^{n+1}M)\]
given by
\[ \bar D \phi(v)=\extd \phi(v)+\omega(v\o)\phi(v\sw \infty).\]
In particular, $T=\bar D\theta$.
\label{prop.dbar}
\end{prop}
\proof We first show that the map $\bar D$ is well-defined. We will use
the following notation for the connection one-form $\omega(c)
= \omega(c)\su 1\tens
\omega(c)\su 2$ (summation understood), for all $c\in C$. Take any
$\phi \in {\rm Hom}_0 (V,P\Omega^n M)$, $v\in V$ and compute
\begin{eqnarray*}
(\id\tens \Delta_{\Omega P})\bar D\phi(v) & = & 1\tens \phi(v)\sw
0\tens \phi(v)\sw 1 + \extd \phi(v)\tens e - 1\tens\phi(v)\tens e \\
&& +\omega(v\sw 1)\su 1\tens \Delta_{\Omega P}(\omega(v\sw 1)\su
2\phi(v\sw \infty))\\
& = & 1\tens \phi(v)\sw
0\tens \phi(v)\sw 1 + \extd \phi(v)\tens e - 1\tens\phi(v)\tens e \\
&& +\omega(v\sw 1)\su 1\tens \omega(v\sw 1)\su 2\sw 0\psi^{n+1}
(\omega(v\sw 1)\su
2\sw 1\tens \phi(v\sw \infty))\\
& = & 1\tens \phi(v)\sw
0\tens \phi(v)\sw 1 + \extd \phi(v)\tens e - 1\tens\phi(v)\tens e \\
&& +1\tens \psi^{n+1}(v\sw 1\tens \phi(v\sw \infty)) - 1\tens
\psi^{n+1}(e\tens \phi(v))\\
&& +\omega(v\sw 1)\psi^{n+1}(v\sw 2 \tens \phi(v\sw \infty))\\
& = & 1\tens \phi(v)\sw
0\tens \phi(v)\sw 1 +\extd \phi(v)\tens e - 1\tens\phi(v)\tens e \\
&& +1\tens \phi(v)\tens e - 1\tens\phi(v)\sw 0\tens \phi(v)\sw 1
+\omega(v\sw 1)\phi(v\sw \infty)\tens e\\
& = & \bar D\phi(v)\tens e,
\end{eqnarray*}
where we used that $\Omega P\in \M_{\Omega P}^C(\psi^\bullet)$ to derive
the second equality, then Proposition~\ref{l.strong.r} to derive the
third one and the fact that $\phi \in {\rm Hom}_0 (V,P\Omega^n M)$ to
obtain the fourth equality. This shows that
$\bar D\phi(v) \in P(\Omega^{n+1} M)$.

Next we need to show that $\bar D\phi$ satisfies (\ref{con.phi.n}).
 We have
\begin{eqnarray*}
\psi^{n+2}(v\sw 1\tens \bar D\phi(v\sw \infty)) & = & \psi^{n+2}(v\sw
1\tens \extd\phi(v\sw\infty)) +\psi^{n+2}(v\sw 1\tens \omega(v\sw
2)\phi(v\sw\infty)) \\
& = & (\extd\tens \id)(\psi^{n+1}(v\sw 1\tens \phi(v\sw\infty))\\
&&+\omega(v\sw 1)\psi^{n+1}(v\sw 2\tens \phi(v\sw\infty)))\\
& = & \extd\phi(v)\tens e+\omega(v\sw 1)\extd\phi(v\sw\infty)\tens e =
\bar D\phi(v)\tens e,
\end{eqnarray*}
where we used the covariance of $\extd$ with respect to $\psi^\bullet$,
the fact that $\Omega P\in \M_{\Omega P}^C(\psi^\bullet)$, and the
covariance property of the connection one-form to derive the second
equality. It is an easy exercise to verify that $T = \bar D\theta$.
\endproof

Here $\bar D$ extends $\bar D$ in Section~4 to higher forms. Next,
again following \cite{Ma:rie}, we introduce left strongly tensorial
forms and a quantum metric. Thus, let $V$ be a right $C$-comodule.
A left strongly tensorial $n$-form is a map $\phi:V\to (\Omega^n M)P$
commuting with the right coaction of $C$, where $C$ coacts on
$(\Omega^nM)P$ by $\psi^{n+1}(e\tens w)$.

\begin{prop} Left strongly tensorial forms $\Hom^C(V,(\Omega^nM)P)$ and
$\Hom_M(\bar E,\Omega^nM)$
are isomorphic as left $M$-modules if $P$ is faithfully flat as a
left $M$-module (cf. \cite{Bou:com} for a comprehensive review of the
concept of faithful flatness).
\end{prop}
\proof This can be shown as \cite[Theorem~5.4]{Brz:mod}. Given $\phi\in
\Hom^C(V,(\Omega^nM)P)$ the corresponding $s_{\phi} \in \Hom_M(\bar
E,\Omega^nM)$ is given by $s_{\phi}(\sum_i v^i\tens
u^i)=\sum_{i} \phi(v^i)u^i$, $\sum_i v^i\tens u^i\in \bar E$. Conversely
given $s\in \Hom_M(\bar E,\Omega^nM)$, the corresponding tensorial form
is given by$\phi_s(v) = s(v\tens 1)$.
\endproof

On the other hand, for $V$ a right $C$-comodule we have the
covariant derivative $D$ extending the $D$ in Section~3 to higher
forms.

\begin{prop} If $\omega$ is a strong connection on $P(M,C,e)$ and
$\psi$ is bijective then there is a covariant derivative
\[ D:\Hom^C(V,(\Omega^nM)P)\to \Hom^C(V,(\Omega^{n+1}M)P)\]
given by
\[ D \phi(v)=\extd \phi(v)+(-1)^{n+1}\phi(v\sw 0)\omega(v\o).\]
\label{prop.d}
\end{prop}
\proof This proposition is a coalgebra bundle version of a similar
statement in \cite{Haj:str} for quantum group principal bundles. 
The proof is similar to the proof of
Proposition~\ref{prop.dbar}. Take any right $C$-covariant  $\phi:V\to
(\Omega^nM)P$ and $v\in V$ and compute
\begin{eqnarray*}
({}_{\Omega P}\Delta\tens \id)D\phi(v) & = & (-1)^{n+1}\phi(v)\sw 1\tens
\phi(v)\sw\infty \tens 1 +e\tens \extd \phi(v) +(-1)^ne\tens\phi(v)\tens
1 \\
&& + (-1)^{n+1} {}_{\Omega P}\Delta(\phi(v\sw 0)\omega(v\sw 1)\su
1)\tens \omega(v\sw 1)\su 2\\
& = & (-1)^{n+1}\phi(v)\sw 1\tens
\phi(v)\sw\infty \tens 1 +e\tens \extd \phi(v) +(-1)^ne\tens\phi(v)\tens
1 \\
&& + (-1)^{n+1} \psi^{-n-1}(\phi(v\sw 0)\tens \omega(v\sw 1)\su
1\sw 1)\omega(v\sw 1)\su1\sw \infty\tens \omega(v\sw 1)\su 2\\
& = & (-1)^{n+1}\phi(v)\sw 1\tens
\phi(v)\sw\infty \tens 1 +e\tens \extd \phi(v) +(-1)^ne\tens\phi(v)\tens
1 \\
&& + (-1)^{n+1} (\psi^{-n-1}(\phi(v\sw 0)\tens v\sw 1)\tens 1
-\psi^{-n-1}(\phi(v)\tens e)\tens 1) \\
&&+ (-1)^{n+1}\psi^{-n-1}(\phi(v\sw 0)\tens v\sw 1)\omega(v\sw 2)
\\
& = & (-1)^{n+1}\phi(v)\sw 1\tens
\phi(v)\sw\infty \tens 1 +e\tens \extd \phi(v) +(-1)^ne\tens\phi(v)\tens
1 \\
&& - (-1)^{n+1}\phi(v)\sw 1\tens
\phi(v)\sw\infty \tens 1 -
(-1)^ne\tens\phi(v)\tens 1\\
&& +(-1)^{n+1}e\tens\phi(v\sw 0)\omega(v\sw 1)\\
& =& e\tens D\phi(v).
\end{eqnarray*}
The third equality follows from Proposition~\ref{l.strong.r}. Thus we
deduce that
$D\phi(v) \in (\Omega^{n+1}M)P$. The proof of the covariance of $D\phi$
is analogous to the corresponding part of the proof of
Proposition~\ref{prop.dbar}. \endproof

Finally, when $V$ is a finite-dimensional left $C$-comodule we can
identify $\bar E$ with $\Hom(V,P)_0$ and
$\Hom^C(V^*,(\Omega^nM)P)$ with $(\Omega^nM)P\Box_C V$ and hence
obtain
\[ (\Omega^nM)P\Box_C V\isom \Hom_M(\Hom_0(V,P),\Omega^nM).\]
We can then define, cf \cite{Ma:rie}, a metric on $M$ as an element
\[ \gamma\in (\Omega^1M)P\Box_C V\]
such that the corresponding map $\Hom_0(V,P)\to \Omega^1M$ is an
isomorphism. In the infinite dimensional case we do not have a
bijection between these spaces,  but we still obtain a map
$\Hom_0(V,P)\to \Omega^1M$ from $\gamma$ and can require it to
be suitably nondegenerate. If $P(M,C,\psi,e)$  and $V$ is a frame 
resolution of
$M$ then we can identify $(\Omega^1M)P\Box_C V$ with $\Omega^2M$,
so that $\gamma$ is a 2-form on $M$.

Following \cite{Ma:rie}, we can also define the {\em cotorsion}
$\Gamma\in\Omega^3M$ of the metric as
\[  \Gamma=(\id\tens s_\theta)(D\Box_C\id)(\gamma).\]
Here, since $\gamma$ is left strongly tensorial (and if $D$
corresponds to a strong connection) then $D\gamma$ is also
left-strongly tensorial when viewed as a map on $V^*$. Hence
$(D\Box_C\id)\gamma\in(\Omega^2M)P\Box_C V$ as required here. In
this context one has the following version of $D$ that does not go
through $V^*$.

\begin{prop} If $\omega$ is a strong connection on $P(M,C,\psi,e)$ and
$\psi$ is bijective then there is a covariant derivative
\[ D:(\Omega^nM)P\Box_C V\to (\Omega^{n+1}M)P\Box_CV\]
given by
\[ D (w\tens v)=\extd w\tens v+(-1)^{n+1}w\omega(v\o)\tens v\sw \infty.\]
\label{prop.d2}
\end{prop}
\proof Dual to the proof of Proposition~\ref{prop.d}.
\endproof

Also provided in \cite{Ma:rie} is a general construction for frame
resolutions on quantum group homogeneous bundles $\pi:P\to H$. We
extend this now in the coalgebra setting $\pi:P\to C$, to
embeddable homogeneous spaces. This more general setting is
definitely needed since it includes, for example, the full family
of quantum 2-spheres \cite{Pod:sph} considered in the next section.
The
following proposition generalises \cite[Proposition~4.3]{Ma:rie}
to include this case.

\begin{prop} A quantum embeddable homogeneous space $M$ of
$P$ corresponding to $\pi:P\to C$ has a coalgebra frame resolution
with $V=M^+$, ${}_V\Delta = (\pi\tens \id)\circ\Delta$ and $\theta
:V\to P(\Omega^1 M)$, $\theta: v\mapsto Sv\sw 1\tens v\sw 2$.
\end{prop}
\proof The canonical entwining structure is $\psi(c\tens h)
= h\sw 1\tens\pi(gh\sw 2)$, where $g\in\pi^{-1}(c)$ (cf.\
\cite[Example~2.5]{BrzMa:coa}). Since $\theta(v)\in
P\tens M$, as $M$ is a left $P$-comodule
algebra, we find
\begin{eqnarray*}
\psi^2(\pi(v\sw
1)\tens\theta(v\sw 2)) & = &Sv\sw 3 \tens \psi(
\pi(v\sw 1 Sv\sw 2)\tens v\sw 4)
 = Sv\sw 1\tens \psi(\pi(1)\tens v\sw 2)\\
& = &Sv\sw 1\tens v\sw 2\tens \pi(1)
=  \theta(v)\tens \pi(1).
\end{eqnarray*}
Since $\tilde\chi(\theta(v)) = (Sv\sw 1)v\sw 2\tens \pi(v\sw 3) =
1\tens\pi(v) = 0$ it follows that  $\theta(v)\in P(\Omega^1 M)$. From
the above calculation we conclude that $\phi\in{\rm
Hom}_0(V,P(\Omega^1 M))$. Now,  consider the map $r: \Omega^1
M\to P\tens M$, $r(\sum_im^i\tens
\tilde{m}^i)=\sum_i m^i\tilde{m}^i\sw 1\tens \tilde{m}^i\sw 2$.
Applying $\id\tens \eps$ to $r$ one immediately finds that ${\rm
Im} r
\subset P\tens V$. Similarly, applying the coaction equalising map
for the cotensor product to $r$ one finds that ${\rm Im} r \subset
P\Box_CV$. Finally using the same argument as in
\cite[Proposition~4.3]{Ma:rie} one
 proves that $r$ is the inverse of $s_\theta: P\Box_CV\to \Omega^1M$,
 $s_\theta :\sum_iu^i\tens v^i \mapsto \sum_iu^iSv^i\sw 1
\tens v^i\sw 2$ as required.
\endproof

\section*{\normalsize\sc\centering 6. Monopole on all quantum
2-spheres}\vspace{-.6\baselineskip}
\setcounter{section}{6}
\setcounter{definition}{0}

Let $SU_q(2)$ be the standard matrix quantum group over the field $k=\C$, 
with generators
$\pmatrix{\alpha&\beta\cr\gamma&\delta}$ and relations
$\alpha\beta=q\beta\alpha$, $\alpha\gamma =q\gamma\alpha$, $\alpha\delta
= \delta\alpha + (q-q^{-1})\beta\gamma$, $\beta\gamma = \gamma\beta$,
$\gamma\delta = q\delta\gamma$, $\alpha\delta - q\beta\gamma =1$. 
Let 
$$
\xi = s(\alpha^2 - q^{-1}\beta^2) +(s^2-1)q^{-1}\alpha\beta, \quad \eta =
s(q\gamma^2 - \delta^2)
+(s^2-1)\gamma\delta,
$$
$$ \zeta = s(q\alpha\gamma - \beta\delta)
+ (s^2-1)q\beta\gamma ,
$$
where $s\in [0,1]$.
We define $C=SU_q(2)/J$ where
$J=\{\xi -s,\eta+s,\zeta\}SU_q(2)$ is a coideal. 
We denote by $\pi$ the canonical projection $SU_q(2)\to C$. As shown in
\cite{Brz:hom}, the fixed point subalgebra
under the coaction of $C$ on $SU_q(2)$  is generated by
$\{1,\xi,\eta,\zeta\}$, and can be
identified with $S_{q,s}^2$, the 2-parameter
quantum sphere in \cite{Pod:sph}. The standard quantum sphere discussed
in \cite{BrzMa:gau}
corresponds to $s=0$. It has been recently proved \cite{MulSch:hom} that 
the coalgebra $C$ is spanned by group-like elements.
We begin by finding such a basis of $C$
explicitly.

\begin{prop} Let 
\[ g^+_n=\pi(\prod_{k=0}^{n-1}
(\alpha + q^ks\beta)), \quad g^-_n=\pi(\prod_{k=0}^{n-1}
(\delta - q^{-k}s\gamma)), \quad n=1,2,\ldots
\]
(all products
increase from left to right). Then $g^\pm_n$ are group-like elements  of
$C$, and
$\{e=\pi(1),g^\pm_n \; | \; n\in {\Bbb N}\}$ is a basis of $C$.
\label{prop.group-like}
\end{prop}
 To prove Proposition~\ref{prop.group-like} we will need the following
\begin{lem}
Let $\ra$ denote the right action of $SU_q(2)$ on $C$, induced by $\pi$.
Then:
\begin{equation}
 sg^+_{n+1} = g^+_n\ra(s\delta+q^{-n}\gamma) = 
g^+_n\ra(s\delta+q^{-n}\beta), 
\label{action+}
\end{equation}
and
\begin{equation}
sg^-_{n+1} = g^-_n\ra(s\alpha-q^{n}\gamma) = 
g^-_n\ra(s\alpha-q^{n}\beta).
\label{action-}
\end{equation}
\label{lemma.action}
\end{lem}
\proof 
Using the commutation rules in $SU_q(2)$ one easily verifies that for
all $s\in \Bbb C$, and $n\in \Bbb N$
\begin{equation}
(\alpha+q^{n-1}s\beta)(s\delta+q^{-n}\gamma) =
(s\delta+q^{-n+1}\gamma)(\alpha+q^ns\beta).
\label{identity+}
\end{equation}

Note that the form of $J=\ker\pi$ implies that for all $x\in SU_q(2)$
\begin{equation}
 \pi((s\delta+\gamma) x)=s\pi((\alpha+s\beta) x),\quad 
s\pi((\delta-s\gamma) x)=\pi((s\alpha-\beta)
x).
\label{pi1}
\end{equation}
This, together with the identity (\ref{identity+}) immediately implies 
that (\ref{action+}) holds for $n=1$. Now, assume that (\ref{action+})
is true for an $n>1$. Then, using the definition of $g^+_n$ as
well as (\ref{identity+}) we have:
\begin{eqnarray*}
g^+_n\ra(s\delta+q^{-n}\gamma) & = & g^+_{n-1}\ra
(\alpha+q^{n-1}s\beta)(s\delta+q^{-n}\gamma)\\
& = &
g^+_{n-1}\ra(s\delta+q^{-n+1}\gamma)(\alpha+q^ns\beta) \\
& = & sg^+_n\ra(\alpha+q^ns\beta) = sg^+_{n+1}.
\end{eqnarray*}
Therefore the first of equalities (\ref{action+}) holds for any $n\in
\Bbb N$.
Since 
\begin{equation}
\pi(\beta x)=\pi(\gamma x), \qquad \forall x\in SU_q(2), 
\label{pi2}
\end{equation}
also the second of equalities (\ref{action+}) holds.

Equalities (\ref{action-}) are proven in an analogous way, by using the
following identity 
\[(s\alpha-q^{n-1}\beta)(\delta-sq^{-n}\gamma) =
(\delta-sq^{-n+1}\gamma)(s\alpha -q^n\beta).\]
\endproof

\noindent {\sl Proof of Proposition~\ref{prop.group-like}.} An easy
calculation which uses (\ref{pi1}) verifies that 
$g^+_1$ is group-like. Assume that $g^+_n$ is group-like
for an $n>1$. Using the definition of $g^+_{n+1}$ and this
inductive assumption we have
\begin{eqnarray*}
\Delta g_{n+1}^+ & = & g^+_n\ra\alpha\tens g^+_n\ra\alpha + g^+_n\ra 
\beta\tens g^+_n\ra\gamma +  q^nsg^+_n\ra\alpha\tens g^+_n\ra\beta +
q^ns g^+_n\ra 
\beta\tens g^+_n\ra\delta \\
& = & g^+_n\ra\alpha\tens g^+_n\ra(\alpha +q^ns\beta)
+q^ng^+_n\ra\beta\tens g^+_n\ra(s\delta +q^{-n}\gamma)\\
& = & g^+_n\ra\alpha\tens g^+_{n+1}
+q^nsg^+_n\ra\beta\tens g^+_{n+1} \qquad \qquad \mbox{\rm 
(Lemma~\ref{lemma.action})}\\
& = & g^+_{n+1}\tens g^+_{n+1}.
\end{eqnarray*}
Thus we conclude that $g^+_n$ is group-like for any $n$. Similarly one
proves that all the $g^-_n$ are group-like. The proof that $\pi(1)$,
$g^\pm_n$ span $C$ is analogous to the proof of
\cite[Proposition~6.1]{Brz:hom}. \endproof

Proposition~\ref{prop.group-like}
 gives an explicit description of the coalgebra bundle. We now
construct a bicovariant splitting of $\pi$ and hence a strong
connection on it.

\begin{prop} The map $i:C\to SU_q(2)$ given by
\[ i(g^+_n)=\prod_{k=0}^{n-1}{\alpha + q^ks(\beta+\gamma) + q^{2k}s^2\delta
\over 1+q^{2k}s^2}, \quad
i(g^-_n)=\prod_{k=0}^{n-1}{\delta - q^{-k}s(\beta+\gamma) +
q^{-2k}s^2\alpha \over 1+q^{-2k}s^2}\]
is bicovariant and splits $\pi$.
\end{prop}
\proof An easy direct calculation which uses (\ref{pi1}), (\ref{pi2}), 
verifies that $\Delta_{SU_q(2)}(i(g^+_1))
=i(g^+_1)\otimes g^+_1$ and ${}_{SU_q(2)}\Delta(i(g^+_1))
=g^+_1\otimes i(g^+_1)$. Now assume that there is $n>1$ such that
$i(g^+_n)$ is bicovariant. Then we have
\begin{eqnarray*}
\Delta_{SU_q(2)}(i(g^+_n)) &=& {1\over
1+q^{2n}s^2}\Delta_{SU_q(2)}(i(g^+_{n+1})(\alpha+sq^n(\beta+\gamma)+
s^2q^{2n}\delta))\\
& = & {1\over 1+ q^{2n}s^2}i(g^+_n)(\alpha
\tens g^+_n\ra(\alpha+q^ns\beta)+q^n\beta\tens
g^+_n\ra(s\delta+q^{-n}\gamma)\\
&& +q^ns\gamma\tens
g^+_n\ra(\alpha+q^ns\beta) + q^{2n}s\delta\tens
g^+_n\ra(s\delta+q^{-n}\gamma)) \\
& = & {1\over
1+q^{2n}s^2}i(g^+_n)(\alpha+q^ns(\beta+\gamma)+q^{2n}s^2\delta)\tens
g^+_{n+1} \qquad \mbox{\rm (Lemma~\ref{lemma.action})}\\
& = & i(g^+_{n+1}) \tens g^+_{n+1}.
\end{eqnarray*}
For the left coaction  we have
\begin{eqnarray*}
{}_{SU_q(2)}\Delta(i(g^+_{n+1})) &=& {1\over
1+q^{2n}s^2}{}_{SU_q(2)}\Delta(i(g^+_n)(\alpha+sq^n(\beta+\gamma)+
s^2q^{2n}\delta))\\
& = & {1\over 1+ q^{2n}s^2}(
g^+_n\ra(\alpha+q^ns\gamma)\tens i(g^+_n)\alpha+q^n
g^+_n\ra(s\delta+q^{-n}\beta)\tens i(g^+_n)\gamma\\
&& +q^ns
g^+_n\ra(\alpha+q^ns\gamma)\tens i(g^+_n)\beta + q^{2n}s
g^+_n\ra(s\delta+q^{-n}\beta)\tens i(g^+_n)\delta) \\
& = & {1\over
1+q^{2n}s^2}g^+_{n+1}\tens i(g^+_n)(\alpha+q^ns(\beta+\gamma)
+q^{2n}s^2\delta)
\qquad \mbox{\rm (Lemma~\ref{lemma.action})}\\
& = & g^+_{n+1} \tens i(g^+_{n+1}).
\end{eqnarray*}
Thus we conclude that $i(g^+_n)$ is bicovariant for all $n\in\Bbb N$.
Similarly one shows that $i(g^-_n)$ is bicovariant. 

The fact that $i$ splits $\pi$ can be proven inductively too, and in the
proof one uses Lemma~\ref{lemma.action} and (\ref{pi1}), (\ref{pi2}). 
\endproof

Consequently, we have a strong connection on $S^2_{q,s}$
defined via the elements $\omega^\pm_n=Si(g_n^\pm)\o\tens
i(g_n^\pm)\t$.

\begin{lem} The elements $\omega^\pm_n$ may be computed iteratively 
from
\[ (1+q^{2n}s^2)\omega^+_{n+1}=(\delta- q^{n+1} s\gamma)
\omega^+_{n}(\alpha+ q^{n}
s\beta)
+(\alpha q^{n}s-
q^{-1}\beta)\omega^+_{n}(q^{n}s\delta+\gamma)\]
\[ (1+q^{-2n}s^2)\omega^-_{n+1}=(q\gamma+q^{-n}s\delta)
\omega^-_{n}(-\beta+q^{-n}s\alpha)
+(\alpha +
q^{-n-1}s\beta)\omega^-_{n}(\delta-q^{-n}s\gamma)\]
and $\omega^\pm_0=1\tens 1$. 
\end{lem}
\proof From $i(g_{n+1}^+)=i(g_{n}^+)(\alpha+sq^n(\beta+\gamma)
+q^{2n}s^2\delta)/(1+q^{2n}s^2)$
as in the proof of Proposition~5.3 and the coproduct and antipode 
$S$ of $SU_q(2)$ one has
\align{(1+q^{2n}s^2)\omega_{n+1}^+\equad &&=S\alpha\omega^+_n\alpha
+S\beta\omega^+_n\gamma+q^ns(
S\beta\omega^+_n\delta+S\alpha\omega^+_n\beta+S\gamma\omega_n^+
\alpha+S\delta\omega^+_n\gamma)
\\
&&\quad +q^{2n}s^2(S\delta\omega_n^+\delta+S\gamma\omega_n^+\beta)\\
&&=\delta\omega_n^+\alpha-q^{-1}\beta\omega_n^+\gamma-q^{n-1}s\beta
\omega_n^+\delta+q^ns\delta\omega_n^+\beta-q^{n+1}s\gamma\omega_n^+
\alpha+q^ns\alpha\omega_n^+\gamma\\
&&\quad+ q^{2n}s^2\alpha\omega_n^+\delta-q^{2n+1}s^2\gamma\omega_n^+
\beta}
which we then factorise as shown. 

The computation for $\omega_{n+1}^-$ is similar. Actually, when 
$s\ne 0$ we may collect the two cases together as
\[ (1+q^{2n}s^{\pm 2})\omega^\pm_{n+1}=(\delta\mp q^{n+1} s^{\pm 1}
\gamma)\omega^\pm_{n}(\alpha\pm q^{n}
s^{\pm 1}\beta)
+(\alpha q^{n}s^{\pm 1}\mp
q^{-1}\beta)\omega^\pm_{n}(q^{n}s^{\pm 1}\delta\pm\gamma).\]
\endproof

For example,
\[ \omega(g_1^+)={1\over 1+s^2}((\delta-qs\gamma)\extd(\alpha+s\beta)+(
\alpha s-q^{-1}\beta)\extd(\gamma+s\delta))=\omega^+_1-1\tens 1.\]
A closed expression for $\omega$ on all $g_n^\pm$ is possible for 
nonuniversal differential calculi where commutation relations exist 
between differential forms and elements
of $S_{q,s}^2$, along the lines of \cite{BrzMa:gau} for the standard 
$q$-monopole. 

Finally, as an example of an associated bundle, let $V=\C$ with the 
right $C$-comodule
structure $\Delta_V(1)=1\tens g_1^+$. Here and in what follows we 
identify linear maps
from $\C$ with their values at $1\in\C$. Then the  space of strongly 
tensorial zero-forms in Proposition~4.8 can be
computed as
\[ \Hom^C(V,P)=\{u\in P|\ \Delta_Ru=u\tens g_1^+\}=\{x(\alpha+s\beta)
+y(\gamma+s\delta)|\ x,y\in S^2_{q,s}\}.\]
The covariant derivative $D:\Hom^C(V,P)\to\Hom^C(V,(\Omega^1M)P)$ can 
be computed as
\begin{eqnarray}
Du&=&\extd u-u\omega(g_1^+)=1\tens u-u\omega_1^+\nonumber \\
& = &1\tens u
-{u\over 1+s^2}(\delta-qs\gamma,\alpha s-q^{-1}\beta)\tens
\pmatrix{\alpha+s\beta\cr \gamma+s\delta} \label{Dmon} 
\end{eqnarray}
from the form of $\omega^+_1$. Here a matrix product (or vector-covector 
contraction) notation is used.

These $\Hom$-spaces correspond to sections of a bundle $\bar E$. From 
another point of view, we may consider $V_L=\C$ with the left 
coaction ${}_V\Delta(1)=g_1^+\tens 1$
and identify the associated bundle $E=P\square_CV_L=\Hom^C(V,P)$ as the same 
space as above. Similarly, we identify $\Omega^1M\tens_M E
=\Hom^C(V,(\Omega^1M)P)$. From this point of view we can consider the above
covariant derivative as a map $\nabla:E\to \Omega^1M\tens_M E$. 

Finally, from the form of $E$ given above it is clear that $E$ is a rank 2 
projective module over $S^2_{q,s}$ along the same lines as the recent result over the standard 
$q$-sphere in \cite{HajMa:pro}. We use the relation 
\begin{equation}
(\delta -qs\gamma)(\alpha+s\beta) +
(s\alpha-q^{-1}\beta)(\gamma+s\delta) = 1+s^2,
\end{equation}
holding in $SU_q(2)$ to verify that 
\[ {\bf p} = {1\over 1+s^2}
\pmatrix{1-\zeta &\xi\cr -\eta & s^2+q^{-2}\zeta} = {1\over 1+s^2}
\pmatrix{\alpha+s\beta
\cr \gamma+s\delta}(\delta - qs\gamma, s\alpha-q^{-1}\beta)\]
obeys ${\bf p}^2={\bf p}$ as an $S_{q,s}^2$-valued $2\times 2$-matrix, 
and that $(S_{q,s}^2)^2{\bf p}=E$
by the identification of $(x,y){\bf p}$ with $u=x(\alpha+s\beta)+y(\gamma+s\delta)$. 
In terms of this, (\ref{Dmon}) becomes
\[ \nabla((x,y){\bf p})=1\tens(x,y){\bf p}-(x,y){\bf p}\tens {\bf p}
=(\extd(x,y)){\bf p}+(x,y)
(\extd {\bf p}){\bf p}=(\extd (x,y){\bf p}){\bf p}\]
so that $\nabla$ is the Grassmannian connection associated to the 
projective module. 
Further details of the projector computation will be presented 
elsewhere. A similar
result holds for general $n$ along the lines for the standard 
$q$-monopole in \cite{HajMa:pro}.

\section*{\normalsize\sc\centering
Acknowledgements}
Research was originally supported by the EPSRC grant GR/K02244 and in 
the case of TB by a Lloyd's of London Tercentenary Foundation Fellowship
in the later stages.

\end{document}